\documentclass[lettersize,journal]{IEEEtran}

\usepackage{authblk}
\usepackage{hyperref}
\usepackage{amsmath,amssymb,amsfonts}
\usepackage{bm}
\usepackage{algorithm,algorithmic,bbm}
\usepackage{graphicx}
\usepackage{color}
\usepackage{mathrsfs}
\usepackage{booktabs}
\usepackage{pifont}
\usepackage{makecell}
\usepackage{epstopdf}  

\usepackage{hyperref}
\usepackage{url}

\usepackage{xr}
\usepackage{multirow}
\usepackage{indentfirst}
\usepackage{balance}
\usepackage{cite}
\usepackage{caption,subcaption}
\newtheorem{corollary}{Corollary}

\newtheorem{assumption}{Assumption}

\newtheorem{remark}{Remark}

\newtheorem{lemma}{Lemma}
\newtheorem{theorem}{Theorem}

\newcommand{\bsx}{\boldsymbol{x}}
\newcommand{\bsv}{\boldsymbol{v}}
\newcommand{\bsX}{\boldsymbol{X}}

\newcommand{\sumn}{\sum\limits_{i=1}^n}
\newcommand{\sumjn}{\sum\limits_{j=1}^n}
\newcommand{\sumT}{\sum\limits_{t=1}^T}

\allowdisplaybreaks[4]
\begin{document}

\title{Improved Dynamic Regret of Distributed Online Multiple Frank-Wolfe Convex Optimization}
\author{Wentao~Zhang,
        Yang~Shi,~\IEEEmembership{Fellow,~IEEE,}
        Baoyong~Zhang,~\IEEEmembership{Member,~IEEE,}
        Deming~Yuan,~\IEEEmembership{Member,~IEEE,}

\thanks{\emph{Corresponding author: Baoyong Zhang.}

Wentao Zhang, Baoyong Zhang and Deming Yuan  are with School of Automation,  Nanjing University of Science and Technology,
        Nanjing 210094, Jiangsu, P. R. China (e-mail: iswt.zhang@gmail.com, baoyongzhang@njust.edu.cn, dmyuan1012@gmail.com).

        Yang Shi is with the Department of Mechanical Engineering, University of Victoria, Victoria, BC V8W 2Y2, Canada (e-mail: yshi@uvic.ca).}
}

\maketitle

\begin{abstract}
 In this paper, we consider a distributed online convex optimization problem over a time-varying multi-agent network. The goal of this network is to minimize a global loss function through local computation and communication with neighbors.  To effectively handle the optimization problem with a high-dimensional and structural constraint set,
 we develop a distributed online multiple Frank-Wolfe algorithm to avoid the expensive computational cost of  projection operation. The dynamic regret bounds are established as $\mathcal{O}(T^{1-\gamma}+H_T)$ with the linear oracle number $\mathcal{O} (T^{1+\gamma})$, which depends on the
horizon (total iteration number) $T$, the function variation $H_T$, and the tuning parameter $0<\gamma<1$.  In particular, when the prior knowledge of $H_T$ and $T$ is available, the bound can be enhanced to $\mathcal{O} (1+H_T)$.
Moreover, we illustrate the significant advantages of the multiple iteration technique and reveal a trade-off between dynamic regret bound, computational cost, and communication cost. Finally, the performance of our algorithm is verified and compared through the distributed online ridge regression problems with two constraint sets.
 \end{abstract}


\section{Introduction}
 \IEEEPARstart{O}{nline} Convex Optimization (OCO) \cite{shalev2011online, hazan2016introduction} as a powerful paradigm for learning has recently received widespread attention in some complex scenarios of optimization and machine learning. Under the OCO framework, the decision maker first makes an action at each round. Then the decision maker suffers a loss from the environment or adversary and receives the information about the loss function to update the next action. The main goal of the decision maker is to minimize  the loss accumulated over time.

In the past two decades, the related work of distributed optimization has increased rapidly \cite{nedic2008distributed, nedic2018distributed, yang2019survey, li2022survey, 10035518,liu2020unitary, li2021distributed, xu2016distributed,6311406, hosseini2013online, yuan2022distributedauto, yi2020distributedtsp, 9806334, shahrampour2017distributed,yuan2020distributed} due to  its prominent advantages of a low computational burden, the robustness as opposed to centralized  structure, and its wide applications, such as sensor networks, signal processing, smart grids, machine learning, etc.  Incorporating the OCO framework, \cite{6311406} and \cite{hosseini2013online} early presented the convergence analysis of distributed online optimization  and established an $\mathcal{O} (\sqrt{T})$ regret for convex loss function. After that, a line of distributed algorithms suitable for the OCO framework, such as distributed online mirror descent \cite{shahrampour2017distributed,yuan2020distributed}, distributed online push-sum \cite{wang2020push}, were explored successively.
 This paper considers the following distributed optimization problem under the OCO framework.
\begin{flalign} \label{problem definition}
\begin{array}{cc}
 \min\limits_{\bsx_1, \cdots, \bsx_T} &\sum\limits_{t=1}^T F_t (\boldsymbol{x}_t) \\
 \text{s.t.} & \boldsymbol{x}_t \in \boldsymbol{X}
\end{array}
\end{flalign}
where $T$ is the time horizon, $\bsX \in {\mathbb{R}}^d$ is a convex compact set, the function $F_t (\boldsymbol{x}_t)= \sum_{i=1}^n{f_{i,t}}( \boldsymbol{x}_t )$, and $f_{i,t}$ is convex in $\bsX$.
Usually, static regret and dynamic regret are used as two performance metrics to measure online optimization algorithms. Static regret $\textbf{Regret}_{s}^j (T)$ \cite{zinkevich2003online} shown in (\ref{static regret}) represents the difference between the cumulative loss incurred from the decision sequence $\{ \boldsymbol{x}_{j,t}\}$ of the agent $j$ over time $T$ and the total loss at the optimal benchmark $ \bsx^* $.
\begin{flalign} \label {static regret}
\textbf{Regret}_{s}^j (T) =   \sum\limits_{t=1}^T F_t { (\boldsymbol{x}_{j,t})} -  \sum\limits_{t=1}^T F_t { (\boldsymbol{x}^*)} &
\end{flalign}
where $
\boldsymbol{x}^* \in {\arg\min}_{\boldsymbol{x} \in \boldsymbol{X}}    \ \sum_{t=1}^T F_t (\boldsymbol{x}).$
In contrast, as a  more stringent metric, dynamic regret \cite{zinkevich2003online, besbes2015non}  defined in (\ref{Regret-j}) accurately reflects the quality of decisions in many applications, such as the objective tracking with multiple robots \cite{xu2022online},  due to its varying benchmark sequence $\{\bsx_1^*, \bsx_2^*, \ldots, \bsx_T^* \}$.
\begin{flalign} \label {Regret-j}
\textbf{Regret}_{d}^j (T) =   \sum\limits_{t=1}^T F_t { (\boldsymbol{x}_{j,t})} -  \sum\limits_{t=1}^T F_t { (\boldsymbol{x}_t^*)} &
\end {flalign}
where $
\boldsymbol{x}_t^* \in {\arg\min}_{\boldsymbol{x} \in \boldsymbol{X}}    \ F_t (\boldsymbol{x}).$

 \begin{table*}[t] 
 \label{table_1}
 \renewcommand\arraystretch{1.5}
\begin{center}
\captionsetup{font={small}, justification=raggedright}
  \caption{The comparison of related works on online Frank-Wolfe optimization.}
  \resizebox{\linewidth}{!}{ 
  \label{table_1}
 \begin{tabular}{c|c|c|c|c|c|c}
\hline
 Reference & Loss function&\makecell[c]{Problem type} & \makecell[c]{Performance \\metric (regret) } &\makecell[c]{Linear oracle } & Regret bound&  \makecell[c]{Requiring prior \\ knowledge (type)$^{\dag}$}\\ \hline   
 Zhang \emph{et al.} \cite{pmlr-v70-zhang17g} & Convex  &Distributed & Static   &  $\mathcal{O}(T)$ & $\mathcal{O}(T^{3/4})$ &Yes $(T)$\\\hline
 Wan \emph{et al.} \cite{pmlr-v119-wan20b} & Convex  &Distributed & Static &$\mathcal{O}(\sqrt{T})$ & $\mathcal{O}(T^{3/4})$ & Yes $(T)$\\\hline
 Wan \emph{et al.} \cite{wan2021projection} & Strongly convex &Distributed  & Static &$\mathcal{O}(T^{1/3})$& $\mathcal{O} (T^{2/3}\log T)$ &Yes $(T)$\\\hline
 Th$\acute{\breve{a}}$ng \emph{et al.} \cite{thuang2022stochastic} & Convex &Distributed & Static  &$\mathcal{O} (T^{3/2})$& $\mathcal{O} (\sqrt{T})$ &Yes $(T)$\\\hline
 \multirow{2}{*}{ Kalhan \emph{et al.} \cite{kalhan2021dynamic}} & \multirow{2}{*}{Convex} & \multirow{2}{*}{Centralized} &\multirow{2}{*}{Dynamic } & $\mathcal{O}(T)$ & $\mathcal{O}(\sqrt{T}(1+H_T+\sqrt{D_T}))$ & \multirow{2}{*}{Yes $(T)$}\\
                          &  & &  & $\mathcal{O}(T^{3/2})$ & $\mathcal{O}(1+H_T+\sqrt{T})$ & \\\hline
 \multirow{2}{*}{ Wan \emph{et al.} \cite{Wan_Xue_Zhang_2021}} & Convex &Centralized & Dynamic &$\mathcal{O}(T\log_2 T)$ & $\mathcal{O} (\max\{\sqrt{T}, T^{2/3}H_T^{1/3}\})$ & \multirow{2}{*}{Yes $(T)$} \\
 &Strongly convex &Centralized & Dynamic  &$\mathcal{O}(T\log_2 T)$& $\mathcal{O} (\max\{\sqrt{TH_T\log T}, \log T\})$ &\\\hline
 Zhang \emph{et al.} \cite{zhang2023dynamic}& Convex &Distributed & Dynamic &$\mathcal{O} (T)$ & $\mathcal{O}(\sqrt{T(1+H_T)}+D_T)$ &Yes $(H_T$ and $T)$ \\\hline
 \textbf{This work } & Convex &Distributed &  Dynamic  &  $\mathcal{O} (T^{1+\gamma})$& $ \begin{array}{cc}
    \mathcal{O}(T^{1-\gamma}+H_T),   0<\gamma<1 \\
    \mathcal{O}(1+H_T),     1- \log_T H_T \leq \gamma \leq 1  \\
\end{array} $ &$ \begin{array}{c}
     \text{No} \\
    \text{Yes} \ (H_T \ \text{and} \  T) \\
\end{array}$ \\\hline
\end{tabular}
}
\begin{flushleft}
$\dag$ Note: In order to obtain regret results, some input parameters of algorithms, such as step size, may require certain prior knowledge.
\end{flushleft}
\end{center}
\vspace{-1.5em}
\end{table*}

Usually, the bound of dynamic regret relies on certain regularity of the optimization problem  and is not sublinear  unless the variation budget satisfies  sublinearly growth in $T$ \cite{besbes2015non}. With that in mind, \emph{the function variation} $H_T$  is defined as
\begin{flalign}
H_T&= \sum_{t=1}^{T-1} \max\limits_{i \in \mathcal{V}} \max\limits_{\bsx \in \boldsymbol{X}} |  f_{i,t+1}(\bsx)- f_{i,t}(\bsx)|  \label {H_T}.
\end{flalign}

This requirement that $H_T$ satisfies sublinear growth in $T$ implies that with the algorithm running, the variability of the function difference at the same $\bsx$ decreases over time. Otherwise, the irregular variability of function difference will lead to $H_T$ being at least of order $T$ because there is always a constant $a>0$ such that $$\max\limits_{i \in \mathcal{V}} \max\limits_{\bsx \in \boldsymbol{X}} |  f_{i,t+1}(\bsx)- f_{i,t}(\bsx)|>a$$ holds. Further, according to \cite{besbes2015non}, this must incur  a dynamic regret of order $T$ under any admissible strategy, and the developed algorithm can not achieve the long-run-average  optimality.
Thus, this requirement is appropriate and required for ensuring that the considered dynamic regret is sub-linear, which can be verified again from the main results given in the following sections.

In such an algorithm design of distributed online optimization problem defined in (\ref{problem definition}), projection operations are usually regarded as  a fundamental method to address the constraints,
such as distributed online gradient descent in \cite{6311406}, distributed online mirror descent in \cite{shahrampour2017distributed}, distributed online dual averaging in \cite{hosseini2013online}.
Generally speaking, the calculation amount brought by the projection operators is equivalent to solving a convex quadratic problem \cite{hazan2012projection}. However, for some optimization scenarios with  a high-dimensional
and structural constraint sets, such as   semidefinite programs in  \cite{hazan2008sparse}, multiclass classification in \cite{pmlr-v70-zhang17g}, image reconstruction \cite{harchaoui2015conditional}, matrix completion in \cite{hazan2012projection,7883821}, matching
pursuit \cite{locatello2017unified}, and optimal control in \cite{wu1983conditional},
the projection step incurs an expensive computational cost. In contrast, Frank-Wolfe (FW) algorithm, also called conditional gradient or projection-free algorithm,  has a significant advantage in saving computational cost \cite{7883821, pmlr-v70-zhang17g} due to the \emph{linear oracle} ${\arg\min}_{\boldsymbol{x} \in \boldsymbol{X}} \langle \boldsymbol{x}, \boldsymbol{b} \rangle$ rather than a projection operator, where $\boldsymbol{b}$ is a known vector.
\subsection{Related Work}
In \cite{hazan2012projection}, Hazan and Kale carried out an early study for centralized online FW optimization, in which the bound of static regret was analyzed. After that,  more and more research works on FW algorithm under the centralized and distributed OCO frameworks  were carried out. In the following, we present a brief review of the works most closely related to this paper and make a comparison in detail in Table \ref{table_1}.

In \cite{pmlr-v70-zhang17g}, Zhang \emph{et al.}  earlier developed an FW algorithm paradigm satisfying distributed online convex optimization (DOCO) framework and obtained the upper bound in $\mathcal{O} (T^{3/4})$ of static regret. Then, \cite{pmlr-v119-wan20b} and \cite{wan2021projection} considered  several improved variants under a low communication frequency and obtained the static regret bounds in $\mathcal{O} (T^{3/4})$ and $\mathcal{O} (T^{2/3}\log T)$ under the conditions of the convex and strongly convex  loss function, respectively.  Th$\acute{\breve{a}}$ng \emph{et al.} \cite{thuang2022stochastic} analyzed two algorithm versions of exact and stochastic gradient  under smooth loss function and showed the static regret bound in $\mathcal{O}(\sqrt{T})$.

For the  more stringent metric: dynamic regret, there is little research work for online FW algorithms, especially in distributed scenarios. Wan \emph{et al.} \cite{Wan_Xue_Zhang_2021} proposed a centralized online FW algorithm combined with a restarting strategy and analyzed dynamic regret  bound for convex and  strongly convex functions, respectively. In \cite{kalhan2021dynamic}, Kalhan \emph{et al.} analyzed the dynamic regret of an online FW algorithm under the condition of a smooth loss function and further improved the dynamic regret bounds by using multiple iterations. For the distributed scenarios, Zhang \emph{et al.} \cite{zhang2023dynamic} developed an online distributed FW optimization algorithm by combing gradient tracking technique and established the dynamic regret bound in $\mathcal{O}(\sqrt{T(1+H_T)}+D_T)$ with a non-adaptive step size, where $D_T$ represent the \emph{gradient variation}.
However, this optimal bound  depends on two variations of the optimization problem and requires the step size that contains a hard-to-get prior knowledge of $H_T$, thus resulting in  difficulties in accurate tuning of the step size parameter  in practice.
\subsection{Motivations and Challenges}
The above analysis and the existing results in Table \ref{table_1} expose the limitations and deficiencies of distributed online FW algorithm in dynamic regret analysis. Is there a method for the algorithm to both remove the above limitations and also establish a tighter regret bound?

  \emph{Multiple iterations}, a method that seeks a higher quality decision by continuously exploiting the information of the function at time $t$, solves this question. In \cite{zhang2017improved}, Zhang \emph{et al.} considered this multiple iterations method for a centralized online gradient descent algorithm and validated that it  could significantly enhance the dynamic regret bound. Eshraghi and Liang \cite{eshraghi2022dynamic} applied this method to the centralized online mirror descent algorithm. However, both research works require  the condition of  a strongly convex loss function.
  Inspired by \cite{kalhan2021dynamic, zhang2023dynamic, zhang2017improved}, in this paper, we aim to investigate whether the technique of multiple iterations  at round $t$ can remove the dependence of step size on prior knowledge and improve the dynamic regret bound of distributed online FW optimization algorithms under the convex loss function.

  Such an idea  naturally brings three key challenges in algorithm design and technical analysis.
  \begin{itemize}
    \item \cite{zhang2017improved} \cite{eshraghi2022dynamic} explored the multiple iterations method in centralized settings. However, in a distributed setting, it is uncertain how the multiple iterations method is combined with the algorithm.
    \item Since the multiple iterations method introduces  additional inner loops, the original technical analysis will no longer be fully suitable. In particular,  the consistency error of the algorithm and the convergence relationship between the inner and outer loops will become very difficult to analyze and prove.
    \item It is tricky to properly choose the parameters of the multiple iteration method to guarantee a great convergence performance.
  \end{itemize}

  \subsection{Contributions}
  The main contributions of this work are stated as follows.

        (i) Incorporating  a technique of multiple iterations  at each round $t$, we propose a distributed online multiple Frank-Wolfe (DOMFW)  algorithm over a time-varying network topology, which efficiently addresses the optimization problem with a  structural and high-dimensional constraint set. Moreover, based on the method of gradient tracking, the global gradient estimation rather than the gradient of the agent itself is employed to update the next decision.

       (ii) We illustrate that the multiple iteration technique can enhance the dynamic regret bound of the FW algorithm in distributed scenarios. For two inner iterations parameters settings, both dynamic regret bounds $\mathcal{O}(T^{1-\gamma}+H_T)$ with the linear oracle number $\mathcal{O}(T^{1+\gamma})$ are established, where $0<\gamma<1$. Moreover, compared with the existing results in \cite{kalhan2021dynamic,Wan_Xue_Zhang_2021, zhang2023dynamic}, this obtained bound  does not require a step size dependent on the prior knowledge of $H_T,T$ and can become tighter by tuning the parameter $\gamma$.

        (iii) With the prior knowledge of $H_T$ and $T$, the optimal bound $\mathcal{O} (1+H_T)$ is obtained, which is the same as the regret level in \cite{wan2023improved}, where the latter additionally requires that the loss function is strongly convex and the optima $\bsx_t^* \in \bsX$ satisfy $\nabla F (\bsx_t^*)=\boldsymbol{0}$. Moreover, we reveal   a trade-off between dynamic regret bound, computational cost, and communication cost. Finally, the performance of our algorithm is validated and compared by simulating the distributed online ridge regression problems with two constraint sets.
 \subsection{Organization and Notations}
The remaining of the paper is structured as follows. The optimization problem, necessary assumptions, and algorithm design are presented in Section II. The convergence results and some discussions are analyzed and established in Section III. Sections IV and V show the simulation examples and conclusion, respectively.

\textbf{Notation}: ${\mathbb{R}}^n$ represents the $n$-dimensional Euclidean space.  $\mathbb{Z}$ and $\mathbb{Z}_+$ represents the integers and positive integers set, respectively.  The Euclidean  norm  of a vector $\boldsymbol{z}$ is denoted as $\|\boldsymbol{z}\|$. $\boldsymbol{1}$ represents a vector whose elements are all equal to $1$. $\boldsymbol{g}^T$  stands for the transpose of a vector $\boldsymbol{g}$. $[K]$ and $[\bsv]_{i}$ denote the set $\{1,\ldots, K\}$ and the $i$-th element of vector $\bsv$, respectively. $[A_t]_{ij}$ stands for the element in the $i$-th row and $j$-th column of matrix $A_t$. For two scalar sequences $\{b_i, i\in \mathbb{Z}_+\}$ and $\{c_i>0, i\in \mathbb{Z}_+\}$, $b_i=\mathcal{O} (c_i)$, $b_i=o (c_i)$ and $b_i=\omega (c_i)$ represent that there exists the scalars $a_1, a_2$ and $a_3$ such that $b_i\leq a_1 c_i$, $b_i < a_1 c_i$ and $b_i> a_1 c_i$, respectively.

\section{Problem Formulation}
\subsection{The Optimization Problem and Some Assumptions}

 In this work, the network information exchange between $n$ agent (nodes) is carried out in a directed time-varying multi-agent graph $\mathcal{G}_t=\{\mathcal{V},\mathcal{E}_t,A_t \}$, where $\mathcal{V} : = \{1,2,\ldots,n\},  \mathcal{E}_t \subseteq \mathcal{V} \times \mathcal{V},$ and $ A_t \in \mathbb{R}^{n \times n}$ stand for node set, edge set, and weight, respectively. Let $\mathcal{N}_{i}^{\text {in }}(t)=\{j \mid(j, i) \in \mathcal{E}_t\} \cup\{i\}$ denotes as inner neighbor sets of agent $i$, and agent $i$ has permission to receive the information from its neighbor agent in $\mathcal{N}_{i}^{\text {in }}(t)$ through the network communication. Moreover, when $j \in \mathcal{N}_{i}^{\text {in }}(t)$, $[A_{t}]_{ij}>0$ holds, and otherwise $[A_{t}]_{ij}=0$ holds.

Distributed online optimization
problems can be described as an interactive loop with the environment or  the adversary as follows:
\begin{itemize}
  \item [(1)] each agent $i$ first commits a decision $\bsx_{i,t} \in \boldsymbol{X}$ with the information reference of its neighbors  at every round $t$;
  \item [(2)] the loss $f_{i,t}$ of agent $i$ and its gradient information are revealed by the environment or  the adversary;
  \item [(3)] agent $i$ updates the next decision $\bsx_{i,t+1}$ through using the information about loss function $f_{i,t}$.
\end{itemize}

The goal of this paper is to ensure that the dynamic regret of each agent $i$ achieves  sublinear convergence by designing an effective online distributed optimization algorithm, i.e., $\lim_{T \rightarrow \infty} (\textbf{Regret}_d^j (T) / T) =0, \forall j \in \mathcal{V}$.

Throughout the paper, we make the following assumptions.
\begin{assumption}\label{assump: network}
(a)  When $[A_{t}]_{ij}>0, t \in [T]$, there exists a positive scalar $ \zeta $ such that $ [A_{t}]_{ij}>\zeta$.
(b)  The graph $ \mathcal{G}_{t}$ is strongly connected for all  time $t$.
(c) For  all $ t\in [T]$, $A_t$ is double stochastic, i.e. $\sum_{j=1}^{n}  [A_{t}]_{ij}=\sum_{i=1}^{n} [A_{t}]_{ij}=1$.
\end{assumption}

\begin{assumption}\label{assump: ball}
For the constraint set $\boldsymbol{X}$, there exists a finite diameter $M$, such that, for any $\boldsymbol{x},\boldsymbol{z} \in \boldsymbol{X}$,
$$\max_{\boldsymbol{x},\boldsymbol{z} \in \boldsymbol{X}} \| \boldsymbol{x} -\boldsymbol{z} \| \leq M.$$
\end{assumption}
\begin{assumption}\label{assump: lips funon}
 For all $i\in [n]$ and $t\in[T]$, the function $f_{i,t}$ is Lipschitz continuous on constraint set $\bsX$ with a known positive constant $L_X$,  i.e.,
 \begin{flalign}\label{assu lips fun eq}
 |f_{i,t}(\boldsymbol{x})-f_{i,t}(\boldsymbol{z})| \leq L_{X} \|\boldsymbol{x}-\boldsymbol{z}\|, \forall \bsx, \boldsymbol{z} \in \bsX.
 \end{flalign}
\end{assumption}
\begin{assumption}\label{assump: lips grad}
For all $i\in [n]$ and $t\in[T]$, the function $f_{i,t}$  has a Lipschitz gradient on constraint set $\bsX$ with a known positive constant ${G_X}$, i.e.,
\begin{flalign}
f_{i,t}(\boldsymbol{x})- f_{i,t}(\boldsymbol{z}) &  \leq \langle \nabla f_{i,t} (\boldsymbol{z}), \bsx-\boldsymbol{z}\rangle \nonumber\\
&\quad + \frac{G_X}{2} \|\boldsymbol{x}-\boldsymbol{z}\|^2, \ \forall \bsx, \boldsymbol{z} \in \bsX.
\end{flalign}
\end{assumption}

\begin{remark}
 Assumptions \ref{assump: ball}, \ref{assump: lips funon} are standard in distributed and centralized  optimization, and similar settings can be seen in \cite{ pmlr-v70-zhang17g, 7883821,besbes2015non}. It should be pointed out that Assumption \ref{assump: ball} can be directly obtained by the properties of compact set $\bsX$. According to Lemma 2.6 in \cite{shalev2011online}, Assumption \ref{assump: lips funon} implies that the gradient is bounded, i.e., $\|\nabla f_{i,t}(\bsx) \|\leq L_X$.  Assumption \ref{assump: lips grad}  is equivalent to the fact that
\begin{flalign}
\| \nabla f_{i,t} (\boldsymbol{x}) - \nabla f_{i,t} (\boldsymbol{z}) \|\leq {G_X} \| \boldsymbol{x}-\boldsymbol{z}\|, \forall \bsx, \boldsymbol{z} \in \bsX.
\end{flalign}
\end{remark}

\subsection{Algorithm DOMFW}
In this subsection,  we first develop Algorithm DOMFW, whose description is shown in Algorithm \ref{algorithm 2}. Different from common online algorithm frameworks, Algorithm DOMFW performs  multiple iterations at each time $t$. In detail, when the decision $\bsx_{i,t}$ of agent $i$ is given, the sequence $ \{ \bsx_{i,t}^1, \bsx_{i,t}^2, \ldots, \bsx_{i,t}^{K_t}, \bsx_{i,t}^{K_t+1}\}$ is generated, where $K_t$ is the iteration number of inner iteration. The inner iterations execute from $\bsx_{i,t}^1=\bsx_{i,t}$ and end at step $\bsx_{i,t+1}=\bsx_{i,t}^{K+1}$ after  performing the consistency, gradient tracking and Frank-Wolfe steps in the inner loop.

\begin{algorithm}[ht]
	\renewcommand{\algorithmicrequire}{\textbf{Initialize:} }
	\caption{ (DOMFW) Distributed Online Multiple Frank-Wolfe Optimization Algorithm }
	\label{algorithm 2}
	\begin{algorithmic}[1]
		\REQUIRE Initial  variables $\boldsymbol{x}_{i,1} \in \boldsymbol{X}$,   parameter $0 < \alpha_t \leq 1,$   and nondecreasing integer sequence $\{K_t\}$.
		
		\FOR {$t=1,2,\cdots,T$}
        \STATE Set $\bsx_{i,t}^1=\bsx_{i,t}$
            \FOR {$k=1,2,\cdots K_t$}
                \FOR {Each agent $i \in \mathcal{V}$}
             \STATE Agent $i$ receives $\boldsymbol{x}_{j,t}^k$ from $j \in \mathcal{N}_{i}^{\text {in }}(t)$, and  updates

                \begin{center}
                \quad $\boldsymbol{\hat{x}}_{i,t}^k= \sum_{j \in \mathcal{N}_{i}^{in}(t)}{[A_t]_{ij} \boldsymbol{x}_{j,t}^k}.$
                           \setlength{\parskip}{0.4em}
                \end{center}
\STATE  Gradient information $\nabla f_{i,t}  (\boldsymbol{x}) $ is revealed and agent $i$ executes gradient tracking steps:
            \IF{$k=1 $}
            \STATE $\overline{\nabla} f_{i,t}^1 =  \nabla f_{i,t}  (\hat{\boldsymbol{x}}_{i,t}^1)$,
            \ELSE
            \STATE $\overline{\nabla} f_{i,t}^k = \widehat{\nabla} f_{i,t}^{k-1} + \nabla f_{i,t}  (\hat{\boldsymbol{x}}_{i,t}^k) - \nabla f_{i,t}  (\hat{\boldsymbol{x}}_{i,t}^{k-1}).$
            \ENDIF

               \begin{center}
               \quad $\widehat{\nabla} f_{i,t}^k= \sum_{j \in \mathcal{N}_{i}^{in}(t)}{[A_{t}]_{ij}} \overline{\nabla} f_{j,t}^k,$
               \end{center}

		\STATE Frank-Wolfe step: agent $i$ calculates $\boldsymbol{v}_{i,t}^k$ by using a linear oracle and updates $\boldsymbol{x}_{i,t}^{k+1}$.

             \begin{center}
               \quad $\boldsymbol{v}_{i,t}^k = \underset{\boldsymbol{x} \in \boldsymbol{X}}{\arg\min} \left<\boldsymbol{x}, \widehat{\nabla} f_{i,t}^k \right>, $

               \quad $\boldsymbol{x}_{i,t}^{k+1}= \hat{\boldsymbol{x}}_{i,t}^k +\alpha_t (\boldsymbol{v}_{i,t}^k- \hat{\boldsymbol{x}}_{i,t}^k).$
              \end{center}

		\ENDFOR

          Agent $i$ updates the decision: $\bsx_{i,t+1}=\bsx_{i,t}^{K_t+1}$ when $k=K_t$ holds.

		\ENDFOR
        \ENDFOR
	\end{algorithmic}
\end{algorithm}

At time $t$, the new decision $\bsx_{i,t+1}$ is updated by exploiting the information of the loss function $F_t (\bsx)$ multiple times, which is similar to the process of distributed off-line optimization only for the time (round) $t$. On the one hand, when $K_t$ is set larger, $\bsx_{i,t+1}=\bsx_{i,t}^{K_t+1}$ is actually closer to the optimum $\bsx_{t}^* \in {\arg\min}_{\boldsymbol{x} \in \boldsymbol{X}}    \ F_t (\boldsymbol{x})$. On the other hand, (\ref{Regret-j}) can be converted to the following inequality.
\begin{flalign} \label{algor principle}
&\textbf{Regret}_{d}^j (T) \nonumber \\
 &=   \sum_{t=2}^T [F_t (\bsx_{j,t}) -  F_t  (\bsx_{t-1}^*) +F_t (\bsx_{t-1}^*) -  F_t  (\bsx_{t}^*)] \nonumber \\
 &\quad + F_1 (\bsx_{j,1}) -  F_1  (\bsx_{1}^*)\nonumber \\
 &=   \sum_{t=1}^{T-1} [F_{t+1} (\bsx_{j,t+1}) -  F_{t+1}  (\bsx_{t}^*)] +\sum_{t=1}^{T-1} [F_{t+1} (\bsx_{t}^*) -  F_t  (\bsx_{t}^*)] \nonumber \\
 &\quad + F_1 (\bsx_{j,1}) -  F_T  (\bsx_{T}^*)\nonumber \\
 &\leq   n L_X \sum_{t=1}^{T-1}  \| \bsx_{j,t+1} - \bsx_{t}^*\| + 2 n H_T +nL_X M
\end{flalign}
where in the last inequality we use (\ref{assu lips fun eq}) and the fact $F_1 (\bsx_{j,1}) -  F_T  (\bsx_{T}^*)= F_1 (\bsx_{j,1})-F_1 (\bsx_T^*)+
\sum_{t=1}^{T-1} F_t (\bsx_T^*)-F_{t+1} (\bsx_{T}^*)\leq nM+H_T$.

It is not hard to note that from the perspective of (\ref{algor principle}),  the dynamic regret is related to the term $\| \bsx_{j,t+1} - \bsx_{t}^*\|$, i.e., the
nearness between $\bsx_{j,t+1}$ and $ \bsx_{t}^*$. Thus, from the view of the principle analysis of the multiple iterations method, it is possible to establish a tighter bound than the existing results in Table \ref{table_1} and  remove  the dependence of step size on prior knowledge.
\section{The Convergence Analysis}
In this section, we analyze and establish in detail the upper bound of dynamic regret for Algorithm \ref{algorithm 2}. Based on this general result, the effect of the choice of algorithm parameters on convergence is discussed in Corollaries \ref{corollary 2} and \ref{corollary 3}. After that, we further show the advantages of our developed algorithm through comparing with some existing results and reveal a trade-off between between obtaining high-quality decision and saving resources. Before that, some preliminaries are first given.
\subsection{Preliminaries}
For convenience, we introduce some notations. Firstly, for the square matrix $A_t$ and the positive integer $K_t$, we let $A_t^{K_t}$ denote the $K_t$-th power of $A_t$, that is
$$
A_t^{K_t}\triangleq \overbrace{A_t\times A_t \times\cdots \times A_t}^{K_t}.
$$
It should be noted that $A_t^0=I$, in which $I$ denote the identity matrix.

Next, for $t\geq s \geq 1$, we introduce $\Phi^K(t,s)$ to denote the following transition matrix:
\begin{flalign}
\Phi(t, s)&\triangleq A_{t} A_{t-1} \ldots A_{s},\label{phi} \\
\Phi^K(t,s)&\triangleq A_t^{K_t}A_{t-1}^{K_{t-1}}\cdots A_{s}^{K_{s}}. \label{phi-K}
\end{flalign}
In addition, we set
$$
\Phi^K(t,t+1)=I.
$$
It is worth noting that $\Phi^K(t,s)$ defined in (\ref{phi-K}) is the product of $\sum_{p=s}^tK_p$ matrices, and each matrix  involved in the product satisfies Assumption \ref{assump: network}. By observing this fact and recalling the convergence property of the transition matrix $\Phi(t,s)$
\footnote{\textbf{Lemma} \cite[Corollary 1]{nedic2008distributed}
Let Assumption \ref{assump: network}  hold. Then, for all $i, j \in \mathcal{V}$, we have $ \left|[\Phi(t, s)]_{i j}-\frac{1}{n}\right| \leq \Gamma_1 \sigma_1^{(t-s)}
$  where $\sigma_1=(1-\zeta / 4 n^{2})$ and $\Gamma_1=(1-\zeta / 4 n^{2})^{-1 }$.
},  we obtain the following condition:
 \begin{flalign} \label{phi-K-b2}
\left|[\Phi^K(t, s)]_{i j}-\frac{1}{n}\right| \leq \Gamma_1 \sigma_1^{\sum_{p=s}^tK_p-1}
\end{flalign}
where $\sigma_1= 1-\frac{\zeta}{4 n^{2}}  $ and $\Gamma_1=\left( 1-\frac{\zeta}{4 n^{2}}\right)^{ -1}$.

Similarly, note that $\Phi^K(t, s+1)A_{s}^{{K_{s}}-l}$ is the product of  $\sum_{p=s}^tK_p-l$ matrices satisfying Assumption \ref{assump: network}, where $1\leq s\leq t$ and   $1\leq l\leq K_s-1$. Thus, we have
 \begin{flalign} \label{phi-K-b3}
\left|[\Phi^K(t, s+1)A_{s}^{{K_{s}}-l}]_{i j} -\frac{1}{n}\right| \leq \Gamma_1 \sigma_1^{\sum_{p=s}^tK_p-l-1}.
\end{flalign}
The conditions in (\ref{phi-K-b2}) and (\ref{phi-K-b3}) will be used in the proof of the following equalities.

 Moreover, to facilitate the proof and analysis, we denote the symbols $\bsx_{avg,t}^k$, $\bsx_{avg,t}$, $\boldsymbol{v}_{avg,t}^k$ and $\delta_{i,t}^k$ as the running average of $\bsx_{i,t}^k$, the running average  of $\bsx_{i,t}$,  the running average  of $\boldsymbol{v}_{i,t}^k$,  and gradient difference of agent $i$, respectively.
\begin{gather}
\left\{\begin{array}{ccl}
\boldsymbol{x}_{avg,t}^k&=&\frac{1}{n} \sum\limits_{i=1}^n \bsx_{i,t}^k,\\
\boldsymbol{x}_{avg,t}&=&\frac{1}{n} \sum\limits_{i=1}^n \boldsymbol{x}_{i,t},\\
\boldsymbol{v}_{avg,t}^k&=&\frac{1}{n} \sum\limits_{i=1}^n \boldsymbol{v}_{i,t}^k,\\
\delta_{i,t}^k&=& \nabla f_{i,t}  \left(\hat{\boldsymbol{x}}_{i,t}^{k+1}\right) - \nabla f_{i,t}  \left(\hat{\boldsymbol{x}}_{i,t}^k\right).
\end{array}\right.
\end{gather}
 Along the above equalities, we can further obtain by combining the relations $\bsx_{i,t}= \bsx_{i,t}^1$ and  $\bsx_{i,t+1}=\bsx_{i,t}^{K_t+1}$ from Algorithm \ref{algorithm 2} that
\begin{flalign}
\bsx_{avg,t}^{1}&=\frac{1}{n} \sumn \bsx_{i,t}^{1}=\bsx_{avg,t},\nonumber \\
\bsx_{avg,t}^{K_t+1}&=\frac{1}{n} \sumn \bsx_{i,t}^{K_t+1}=\bsx_{avg,t+1}.
\end{flalign}
\subsection{Main Convergence Results}
Along with the aforementioned basic conditions, some crucial lemmas and the general convergence results of Algorithm \ref{algorithm 2} are established in this subsection. Lemma \ref{consistency2} and Lemma \ref{grad diffience2} give the upper bound of the consistency error of the state variable and the upper bound of the tracking error of the estimated gradient, respectively. Different from the analysis of general distributed online optimization \cite{zhang2023dynamic,shahrampour2017distributed}, the convergence properties of the two new transition matrices described in (\ref{phi-K-b2}) and (\ref{phi-K-b3}) play a key role in both technical analyses.

\begin{lemma} (Consistency error)\label{consistency2}
Suppose Assumptions \ref{assump: network} and \ref{assump: ball} hold. Let $\{ \bsx_{i,t}\}$ be the decision sequence generated by Algorithm \ref{algorithm 2}. Then, we have for any  $T\geq 2$, $K_t\geq 2$ that
\begin{flalign}\label{consistency2-con1}
 &\quad \sumT \sumn \| {\bsx}_{i,t}- \bsx_{avg,t}\| \nonumber\\
& \leq \frac{n\Gamma_1}{1-{\sigma_1}^{K_1}}  \sum_{i=1}^n  \|  {\bsx}_{i,1}\| + \sumn \| {\bsx}_{i,1}- \bsx_{avg,1}\| \nonumber\\
 &\quad + \left(\frac{n^2M \Gamma_1 }{{\sigma_1}\left(1-{\sigma_1}\right)\left(1-{\sigma_1}^{K_1}\right) }+2nM \right)\sumT  \alpha_{t}.
\end{flalign}
\end{lemma}
\noindent{\em Proof:} See Appendix A.

\begin{lemma} (Tracking error) \label{grad diffience2}
Suppose Assumptions \ref{assump: network} and \ref{assump: lips grad} hold. Let $\{ \widehat{\nabla}^k f_{i,t}\}$ and $\{ {\nabla} f_{i,t} (\hat{\boldsymbol{x}}_{i,t}^k) \}$ be the sequences generated by Algorithm \ref{algorithm 2}. Then, we have for any $T\geq 2$, $K_t\geq 2$ that
\begin{flalign}  \label{grad diffience2-con1}
  &\quad\sumT  \sum_{k=1}^{K_t} \sumn \alpha_t \left\|  \widehat{\nabla} f_{i,t}^k - \frac{1}{n}\nabla F_t ({\bsx}_{avg,t}^k)\right\| \nonumber \\
  &\leq D_1 \sum_{t=1}^T \alpha_t  +D_2 \sumT \alpha_t^2  {K_t}
\end{flalign}
where
$
D_1=(\frac{ 2n \Gamma_1 G_X}{1-{\sigma_1}}+G_X) (nM+ \frac{n\Gamma_1 }{1-{\sigma_1}} \sum_{j=1}^n (\| \boldsymbol{x}_{j,1}\|+M) ) +  \frac{ n^2 \Gamma_1 L_X}{1-{\sigma_1}},
D_2=\frac{2n^2\Gamma_1 G_X M}{1-{\sigma_1}}(\frac{n\Gamma_1}{1-{\sigma_1}}+3)+2nMG_X.
$
\end{lemma}
\noindent{\em Proof:} See Appendix B.

\begin{lemma} (Convergence of  the inner loop) \label{lem F_k}
Suppose Assumptions \ref{assump: ball} and \ref{assump: lips grad} hold. Let $\{ \bsx_{i,t}\}$ be the decision sequence generated by Algorithm \ref{algorithm 2}. Then, we have at time $t$ that
\begin{flalign} \label{Proof lem F_k- 3}
&F_{t} \left(\boldsymbol{x}_{avg,t}^{k+1}\right) - F_{t}(x_{t}^*) \leq  (1-\alpha_t)^k  \left[ F_{t}  \left(\boldsymbol{x}_{avg,t}^1 \right)-F_{t}(x_{t}^*)\right] \nonumber\\
&\quad + 2\alpha_t M \sum_{l=1}^{k}   \sumn \left\| \frac{1}{n}\nabla F_{t} \left(\boldsymbol{x}_{avg,t}^{l}\right)- \widehat{\nabla} f_{i,t}^{l} \right\|    + \frac{n\alpha_t G_XM^2}{2}.
\end{flalign}
\end{lemma}
\noindent{\em Proof:} See Appendix C.

 In Lemma \ref{lem F_k}, the dependence of algorithm convergence at time $t$ on the iterations number $k$ of the inner loop is obtained. It should be pointed out that Lemma \ref{lem F_k} plays a key role in linking the convergence of the inner and outer loops of the algorithm in the overall regret analysis.

 Based on this, through using the definition of $H_T$ and combing Lemmas \ref{consistency2}-\ref{lem F_k} and the key inequality about the term $\left(1-  \alpha_t\right)^{ K_t}$, we establish the dynamic regret bound of Algorithm \ref{algorithm 2} in Theorem \ref{theorem 2}.

\begin{theorem}\label{theorem 2}
Suppose Assumptions \ref{assump: network}-\ref{assump: lips grad} hold. Let $\{ \bsx_{i,t}\}$ be the decision sequence generated by Algorithm \ref{algorithm 2} with  the parameter $\alpha_t={1}/{(\rho K_t)}$, where $\rho \geq 1$ is a constant. Then, for  $T\geq 2$, $K_t\geq 2$ and $j \in \mathcal{V}$,  we obtain that
 \begin{flalign} \label{theorem 2-equ}
& \mathbf{Regret}_d^j(T)\leq E_1+ E_2 H_T+     E_3 \sum_{t=1}^T \frac{1}{K_t}
\end{flalign}
where\\
$
E_1=\frac{n^2 L_{X}\Gamma_1}{1-{\sigma_1}^{K_1}}  \sum_{i=1}^n  \|  {\bsx}_{i,1}\| + n L_{X} \sum_{i=1}^n \| {\bsx}_{i,1}- \bsx_{avg,1}\|+ n L_X M (1- \mathrm{e} ^{-\frac{1}{\rho}})^{-1},  \\
E_2=2n(1- \mathrm{e} ^{-\frac{1}{\rho}})^{-1},\\
E_3= 2M(\frac{1}{\rho}D_1+\frac{1}{\rho^2}D_2)(1- \mathrm{e} ^{-\frac{1}{\rho}})^{-1}+\frac{ n G_X M^2}{2 \rho}(1- \mathrm{e} ^{-\frac{1}{\rho}})^{-1}
\quad+ \frac{n^2 L_{X} M}{\rho}(\frac{n \Gamma_1 }{{\sigma_1}(1-{\sigma_1})(1-{\sigma_1}^{K_1}) }+2 ).
$
\end{theorem}

\noindent{\em Proof:} Setting $k=K_t$ in Lemma \ref{lem F_k} and using the above relation,  we get the following inequality:
\begin{flalign} \label{Proof theorm2-1}
& \quad F_{t} \left(\boldsymbol{x}_{avg,t+1}\right) - F_{t}(x_{t}^*)\nonumber \\
 &\leq  (1-\alpha_t)^{K_t}  \left[ F_{t}  \left(\boldsymbol{x}_{avg,t} \right)-F_{t}(x_{t}^*)\right] \nonumber \\
 &\quad + 2\alpha_t M \sum_{l=1}^{K_t}   \sumn \left\| \frac{1}{n}\nabla F_{t} \left(\boldsymbol{x}_{avg,t}^{l}\right)- \widehat{\nabla} f_{i,t}^{l} \right\|    + \frac{n\alpha_t G_XM^2}{2}.
\end{flalign}
The term on the left-hand side of the above inequality has the following relation:
\begin{flalign*}
 &\quad F_{t} (\boldsymbol{x}_{avg,t+1}) -F_{t}(x_{t}^*)\nonumber \\
&= F_{t+1} (\boldsymbol{x}_{avg,t+1}) - F_{t+1} (\boldsymbol{x}_{t+1}^*)+F_{t} (\boldsymbol{x}_{avg,t+1}) \nonumber \\
 &\quad- F_{t+1}(\bsx_{avg,t+1}) +F_{t+1}(x_{t+1}^*) -F_{t}(x_{t}^*).
\end{flalign*}
Note that
$$F_{t+1}(\bsx_{avg,t+1})- F_{t} (\boldsymbol{x}_{avg,t+1}) \leq nf_{t+1,sup}. $$
Then,  we have
\begin{flalign} \label{Proof theorm2-2}
& F_{t+1} (\boldsymbol{x}_{avg,t+1}) - F_{t+1} (\boldsymbol{x}_{t+1}^*) \nonumber \\
   &\leq nf_{t+1,sup} + F_{t}(x_{t}^*)-F_{t+1}(x_{t+1}^*)  \nonumber \\
 &\quad+ (1-\alpha_t)^{K_t}  \left[ F_{t}  \left(\boldsymbol{x}_{avg,t} \right)-F_{t}(x_{t}^*)\right] \nonumber \\
  &\quad + 2\alpha_t M \sum_{l=1}^{K_t}   \sumn \left\| \frac{1}{n}\nabla F_{t} \left(\boldsymbol{x}_{avg,t}^{l}\right)- \widehat{\nabla} f_{i,t}^{l} \right\|    + \frac{n\alpha_t G_XM^2}{2}.
\end{flalign}
Summing from $t=1$ to $T-1$ on both sides of   (\ref{Proof theorm2-2}), we get
\begin{flalign} \label{Proof theorm2-3}
&\sum\limits_{t=1}^{T-1} \left[F_{t+1} (\boldsymbol{x}_{avg,t+1})-F_{t+1}(\boldsymbol{x}_{t+1}^*) -  F_{t} (\boldsymbol{x}_{avg,t}) + F_{t}(\boldsymbol{x}_{t}^*)\right] \nonumber \\
&\leq n H_T +F_{1}(\boldsymbol{x}_{1}^*) -F_{T}(\boldsymbol{x}_{T}^*) + \frac{ n G_X M^2}{2} \sumT \alpha_t \nonumber \\
 &\quad+{ 2 M}\sum\limits_{t=1}^{T-1} \sum_{l=1}^{K_t}   \sumn\alpha_t \left\| \frac{1}{n}\nabla F_{t} \left(\boldsymbol{x}_{avg,t}^{l}\right)- \widehat{\nabla} f_{i,t}^{l} \right\|   \nonumber\\
&\quad -\sum\limits_{t=1}^{T-1}\left[1-(1-\alpha_t)^{K_t}\right] [F_{t} (\boldsymbol{x}_{avg,t})- F_{t}(\boldsymbol{x}_{t}^*)].
\end{flalign}
Note that
\begin{flalign}
&\quad\sum\limits_{t=1}^{T-1} [F_{t+1} (\boldsymbol{x}_{avg,t+1})-F_{t+1}(\boldsymbol{x}_{t+1}^*) - F_{t} (\boldsymbol{x}_{avg,t}) + F_{t}(\boldsymbol{x}_{t}^*)]\nonumber \\
 &=
\left[F_{T} (\boldsymbol{x}_{avg,T})-F_{T}(\boldsymbol{x}_{T}^*)\right]-\left[F_{1} (\boldsymbol{x}_{avg,1})-F_{1}(\boldsymbol{x}_{1}^*)\right].
\end{flalign}
This, together with (\ref{Proof theorm2-3}), implies
\begin{flalign} \label{Proof theorm2-3-b1}
&\sum\limits_{t=1}^{T-1}\left[1-(1-\alpha_t)^{K_t}\right] [F_{t} (\boldsymbol{x}_{avg,t})- F_{t}(\boldsymbol{x}_{t}^*)] \nonumber \\
&\leq n H_T  +{ 2 M}\sum\limits_{t=1}^{T-1} \sum_{l=1}^{K_t}   \sumn\alpha_t \left\| \frac{1}{n}\nabla F_{t} \left(\boldsymbol{x}_{avg,t}^{l}\right)- \widehat{\nabla} f_{i,t}^{l} \right\| \nonumber \\
 &\quad + \frac{ n G_X M^2}{2} \sumT \alpha_t + F_{1} (\boldsymbol{x}_{avg,1}) - F_{T} (\boldsymbol{x}_{avg,T}).
\end{flalign}
Note that $F_{T} (\boldsymbol{x}_{avg,T})- F_{T}(\boldsymbol{x}_{T}^*)\geq 0$. Thus, it follows from (\ref{Proof theorm2-3-b1})    that
\begin{flalign} \label{Proof theorm2-3-c1}
&\sum\limits_{t=1}^{T}\left[1-(1-\alpha_t)^{K_t}\right] [F_{t} (\boldsymbol{x}_{avg,t})- F_{t}(\boldsymbol{x}_{t}^*)] \nonumber \\
&\leq n H_T  +{ 2 M}\sum\limits_{t=1}^{T-1} \sum_{l=1}^{K_t}   \sumn\alpha_t \left\| \frac{1}{n}\nabla F_{t} \left(\boldsymbol{x}_{avg,t}^{l}\right)- \widehat{\nabla} f_{i,t}^{l} \right\| \nonumber \\
 &\quad + \frac{ n G_X M^2}{2} \sumT \alpha_t + F_{1} (\boldsymbol{x}_{avg,1}) - F_{T} (\boldsymbol{x}_{avg,T}) \nonumber\\
&\quad+\left[1-(1-\alpha_t)^{K_T}\right] [F_{T} (\boldsymbol{x}_{avg,T})- F_{T}(\boldsymbol{x}_{T}^*)]\nonumber \\
&\leq n H_T  +{ 2 M}\sum\limits_{t=1}^{T-1} \sum_{l=1}^{K_t}   \sumn\alpha_t \left\| \frac{1}{n}\nabla F_{t} \left(\boldsymbol{x}_{avg,t}^{l}\right)- \widehat{\nabla} f_{i,t}^{l} \right\| \nonumber \\
 &\quad + \frac{ n G_X M^2}{2} \sumT \alpha_t  + F_{1} (\boldsymbol{x}_{avg,1})  - F_{T}(\boldsymbol{x}_{T}^*) \nonumber \\
&\leq 2n H_T +nL_X M + \frac{ n G_X M^2}{2} \sumT \alpha_t \nonumber \\
 &\quad +{ 2 M}\sum\limits_{t=1}^{T-1} \sum_{l=1}^{K_t}   \sumn\alpha_t \left\| \frac{1}{n}\nabla F_{t} \left(\boldsymbol{x}_{avg,t}^{l}\right)- \widehat{\nabla} f_{i,t}^{l} \right\|
\end{flalign}
where the last inequality is established based on the following inequality:
\begin{flalign} \label{proof lem key 6-a5}
&\quad F_{1} (\boldsymbol{x}_{avg,1}) -F_{T}(\boldsymbol{x}_{T}^*) \nonumber \\
&=F_{1} (\boldsymbol{x}_{avg,1})-F_{1}(\boldsymbol{x}_{T}^*)+\sum_{t=1}^{T-1} [F_{t}(\boldsymbol{x}_{T}^*)-F_{t+1}(\boldsymbol{x}_{T}^*)]\nonumber \\
&\leq nL_X \|\boldsymbol{x}_{avg,1}- \boldsymbol{x}_{T}^*\|+\sum_{t=1}^{T-1} \sumn \left|f_{i,t}(\boldsymbol{x}_{T}^*)-f_{i,t+1}(\boldsymbol{x}_{T}^*)\right|\nonumber \\
&\leq nL_X M +n H_T.
\end{flalign}
When $\alpha_t$ is chosen as $1/(\rho K_t)$, in which $\rho\geq 1$,  we have that
$$
 \left(1- \alpha_t\right)^{\rho K_t}= \left(1- \frac{1}{\rho K_t} \right)^{\rho K_t}\leq \mathrm{e} ^{-1}
$$
where $\mathrm{e}$ is the natural constant.  Thus,
$$
 \left(1-  \alpha_t\right)^{ K_t} \leq \mathrm{e} ^{-\frac{1}{\rho}}
$$
which implies that
$$
1- \left(1-  \alpha_t\right)^{ K_t} \geq 1- \mathrm{e} ^{-\frac{1}{\rho}}>0.
$$
By using this fact, it follows from (\ref{Proof theorm2-3-c1}) that
\begin{flalign} \label{Proof theorm2-3-b3}
&\left(1- \mathrm{e} ^{-\frac{1}{\rho}}\right) \sum\limits_{t=1}^{T}[F_{t} (\boldsymbol{x}_{avg,t})- F_{t}(\boldsymbol{x}_{t}^*)] \nonumber \\
&\leq  2n H_T +nL_XM + \frac{ n G_X M^2}{2} \sumT \alpha_t\nonumber \\
 &\quad+{ 2 M}\sum\limits_{t=1}^{T-1} \sum_{l=1}^{K_t}   \sumn\alpha_t \left\| \frac{1}{n}\nabla F_{t} \left(\boldsymbol{x}_{avg,t}^{l}\right)- \widehat{\nabla} f_{i,t}^{l} \right\|.
\end{flalign}

Based on the above analysis, we are now in position to derive the upper-bound of the regret defined in (\ref{Regret-j}).
\begin{flalign}\label{Proof theorm2-5-b1}
 &\quad\textbf{Regret}_d^j(T) \nonumber \\
& = \sum\limits_{t=1}^T \left[F_t {  (\boldsymbol{x}_{j,t})} - F_t (\boldsymbol{x}_{avg,t})\right]
 + \sum\limits_{t=1}^T \left[F_t {  (\boldsymbol{x}_{avg,t})} - F_t { (\boldsymbol{x}_{t}^*)}\right] \nonumber \\
&\leq nL_{X} \sum\limits_{t=1}^T      \|\boldsymbol{x}_{j,t}-\boldsymbol{x}_{avg,t}\|
 + \sum\limits_{t=1}^T \left[F_t {  (\boldsymbol{x}_{avg,t})} - F_t{ (\boldsymbol{x}_{t}^*)}\right]\nonumber \\
&\leq   n L_{X} \sum_{t=1}^T \sum_{i=1}^n   \|\boldsymbol{x}_{i,t}-\boldsymbol{x}_{avg,t}\|  + \sum_{t=1}^T [F_t {  (\boldsymbol{x}_{avg,t})} - F_t{ (\boldsymbol{x}_{t}^*)}]\nonumber \\
&\leq\left(1- \mathrm{e} ^{-\frac{1}{\rho}}\right)^{-1} \left[  2n H_T +n L_X M +\frac{ n G_X M^2}{2} \sumT\alpha_t  \right.\nonumber \\
  &\left.\quad +{ 2 M} \sum\limits_{t=1}^{T-1} \sum_{l=1}^{K_t}   \sumn \alpha_t\left\| \frac{1}{n}\nabla F_{t} \left(\boldsymbol{x}_{avg,t}^{l}\right)- \widehat{\nabla} f_{i,t}^{l} \right\| \right]    \nonumber \\
&  \quad + n L_{X} \sum_{t=1}^T \sum_{i=1}^n   \|\boldsymbol{x}_{i,t}-\boldsymbol{x}_{avg,t}\|.
%
\end{flalign}
Then, by applying Lemmas  \ref{consistency2} and \ref{grad diffience2}, (\ref{theorem 2-equ}) can be obtained.
The proof is complete.
\hfill$\square$

\subsection{Discussions}
Theorem \ref{theorem 2} shows the main results of dynamic regret for Algorithm DOMFW-CO. It is easy to note that the regret bound of Algorithm \ref{algorithm 2} depends on the choices of $\alpha_t$ or sequence $\{K_t \}$. Hence, we have the following corollary by choosing suitable sequence $\{K_t \}$.
\begin{corollary} \label{corollary 2}
Let the conditions in Theorem \ref{theorem 2} hold.
   Then, if $H_T = o(T)$ holds, taking $K_t={1}/{(\rho \alpha_t)}= \lceil \varepsilon_1 t^{ \gamma_1}\rceil +1, \varepsilon_1 >0, 0<\gamma_1<1, t\in\{1,2, \ldots,T\}$, we have
\begin{flalign}\label{corollary 1 equation}
 \mathbf{Regret}_d^j(T) &\leq
    \mathcal{O}\left(T^{1-\gamma_1}+H_T\right),   \nonumber\\
\textbf{N}_{LO}&\leq \mathcal{O}(T^{1+\gamma_1})
\end{flalign}
where $\textbf{N}_{LO}$ denotes the iteration number of linear oracle ${\arg\min}_{\boldsymbol{x} \in \boldsymbol{X}} \langle \boldsymbol{x}, \widehat{\nabla} f_{i,t}^k \rangle$ in Algorithm \ref{algorithm 2}.

 In particular, when $1- \log_T H_T \leq \gamma_1 <1, H_T=\omega(1)$, the upper bounds of dynamic regret and $\textbf{N}_{LO}$ can be establish as $\mathcal{O}(H_T)$ and $\mathcal{O}(T^{2- \log_T H_T})$, respectively.
\end{corollary}
\noindent{\em Proof.} Substituting the conditions in Corollary \ref{corollary 2} into inequality (\ref{theorem 2-equ}), we obtain
\begin{flalign} \label{Cor_prof1}
 \mathbf{Regret}_d^j(T)&\leq E_1+ E_2 H_T+     E_3 \sum_{t=1}^T \frac{1}{\lceil \varepsilon_1 t^{\gamma_1}\rceil +1} \nonumber \\
&\leq E_1+ E_2 H_T+     \frac{E_3}{\varepsilon_1 } \sum_{t=1}^T \frac{1}{t^{\gamma_1}}\nonumber \\
&\leq \mathcal{O}\left(T^{1-\gamma_1}+H_T\right)
\end{flalign}
where the last inequality is obtained by using
\begin{flalign}
&  \sum_{t=1}^T\frac{1}{t^{\gamma_1}} \leq 1+ \sum_{t=2}^T \frac{1}{t^{\gamma_1}} \leq 1+ \int_1^T{{\frac{1}{t^{\gamma_1}} }} dt \leq \frac{1}{1-\gamma_1} T^{1-\gamma_1}.
\end{flalign}
According to Algorithm \ref{algorithm 2}, $\textbf{N}_{LO}$ is equivalent to the sum of the inner loop numbers $ \sum_{t=1}^T K_t$. Thus, we have that
\begin{flalign} \label{Cor_prof3}
%
%
& \sumT K_t \leq \sumT ( \varepsilon_1t^{\gamma_1}+2)\leq \int_1^T ( \varepsilon_1t^{\gamma_1}+2)dt\leq \mathcal{O}(T^{1+\gamma_1}).
\end{flalign}

In particular, through substituting the conditions $1- \log_T H_T \leq \gamma_1 <1, H_T=\omega(1)$ into (\ref{Cor_prof1}) and (\ref{Cor_prof3}), the bounds of dynamic regret and $\textbf{N}_{LO}$ are naturally obtained. The proof is complete.
\hfill$\square$
\begin{corollary} \label{corollary 3}
Let the conditions in Theorem \ref{theorem 2} hold.
   Then, if $H_T =o(T)$ holds, taking $K_t={1}/{(\rho \alpha_t)}= \lceil \varepsilon_2 T^{ \gamma_2}\rceil +1, \varepsilon_2 >0, 0<\gamma_2\leq 1, t\in\{1,2, \ldots,T\}$, we have
\begin{flalign}\label{corollary 3 equation}
 \mathbf{Regret}_d^j(T)
&\leq     \mathcal{O}\left(T^{1-\gamma_2}+H_T\right),\nonumber\\
\textbf{N}_{LO}&\leq \mathcal{O}(T^{1+\gamma_2}).
\end{flalign}

In particular, when $1- \log_T H_T \leq \gamma_2 \leq 1$, the upper bounds of dynamic regret and $\textbf{N}_{LO}$ can be establish as $\mathcal{O}(1+H_T)$ and $\mathcal{O}(T^{2-  \log_T H_T})$, respectively.
\end{corollary}

\noindent{\em Proof.} Similar to the proof of Corollary \ref{corollary 2},  we have that
\begin{flalign}
 \mathbf{Regret}_d^j(T)&\leq E_1+ E_2 H_T+     \frac{E_3}{\varepsilon_2 } {T^{1-\gamma_2}}, \nonumber \\
  \sumT K_t &\leq \mathcal{O}(T^{1+\gamma_2}).
\end{flalign}
The proof is complete.
\hfill$\square$

\begin{remark} \label{corollary remark}
Specially, the regret bounds in Corollaries \ref{corollary 2} and \ref{corollary 3} match the centralized results in \cite{kalhan2021dynamic} and the distributed results in \cite{zhang2023dynamic}, and have more significant advantages than them.
\begin{enumerate}
  \item[i)] Compared with the bound $\mathcal{O}(\sqrt{T}(1+H_T+\sqrt{D_T}))$ of the first algorithm in \cite{kalhan2021dynamic}, our results are less conservative and tighter, and can remove the dependency on  $D_T$, where $D_T$ denotes \emph{gradient variation}. For example, the bound $\mathcal{O}(\sqrt{T}+H_T)$ is better when $\gamma_1=0.5$ holds. The range of $H_T$ in this paper is further enhanced to $\mathcal{O} (T)$ instead of $\mathcal{O} (\sqrt{T})$ if a sublinear regret bound is expected, which effectively expands the application field of optimization problems.
  \item [ii)] In contrast to the bound $\mathcal{O}(1+\sqrt{T}+H_T)$ of the second algorithm in \cite{kalhan2021dynamic},  the regret bounds in Corollaries \ref{corollary 2} and \ref{corollary 3} are more flexible and can achieve a better convergence performance, i.e., the cases when $\gamma_1\in (0.5,1)$ or $\gamma_2\in (0.5,1]$.
      Moreover, the step size in Corollary \ref{corollary 2} does not require the prior knowledge of  $T$.
  \item [iii)] In the distributed work  \cite{zhang2023dynamic}, the establishment of the optimal dynamic regret bound $\mathcal{O} (\sqrt{T(1+H_T)} +D_T)$ depends on $D_T$ and a step size with knowledge of $H_T$, which leads to difficulties in accurate tuning of the step size parameter  in practice. In contrast, our results remove the aforementioned limitations and have better convergence performance than \cite{zhang2023dynamic}.
\end{enumerate}
\end{remark}

\begin{remark} (Optimal bound) Corollaries \ref{corollary 2} and  \ref{corollary 3} reveal that under the prior knowledge of $H_T$ and $T$, the general bounds in (\ref{corollary 1 equation}) and (\ref{corollary 3 equation}) can be improved to the optimal bound $\mathcal{O}(1+H_T)$, which also can achieve a saving in computing and communication resources by   a tunning  of lower inner loop number $K_t$. In particular, this optimal bound $\mathcal{O} (1+H_T)$ is the same as the regret level in \cite{wan2023improved}, where the latter requires that the loss function is strongly convex and the optima $\bsx_t^* \in \bsX$ satisfy $\nabla F (\bsx_t^*)=\boldsymbol{0}$. Therefore, this also reflects that $H_T$ has a large impact on dynamic regret.
Moreover, in practice, the order of $H_T$ over time $T$ may be difficult to determine exactly, and its estimation $\hat{H}_T:=\mathcal{O} (T^{\gamma_3})\geq H_T $ based on factitious experience is undoubtedly a practical substitute, where $\gamma_3$ is a known constant.

\end{remark}
\begin{remark} (Communication number) It is not hard to note that from Algorithm \ref{algorithm 2}, agent $i$ communicates $2K_t$ times with its neighbor time $t$. As the algorithm runs over time $T$,  the level of communication number attaches $\mathcal{O}(T^{1+\gamma_1})$ under the parameter settings of Corollary \ref{corollary 2}, which is a weakness of this paper since more communication resources are used than non-multiple algorithms. At time $t$, the multiple communications between agents is to ensure that the algorithm can obtain detailed neighbor information to output a more reliable decision $\bsx_{i, t+1}$. Thus, from this point of view, it is reasonable.
\end{remark}

\subsection{The Trade-off Between Regret Bound, Computational Cost,  and Communication Cost } \label{trade-off}
From Corollaries \ref{corollary 2} and \ref{corollary 3}, a trade-off between regret bounds and resource consumption can be revealed.
On the one hand, it is easy to find that as the parameter $\gamma_1\in (0,1) $ or $\gamma_2\in (0,1]$ increases, the dynamic regret bound becomes tighter. In other words, when the scale $k$ of decision-making develops much faster than the time-scale $t$ of an online optimized process, a large parameter $\gamma_1$ or $\gamma_2$ is beneficial for Algorithm \ref{algorithm 2} to obtain high-quality decisions.
On the other hand, such a great regret level does not mean that a larger parameter $\gamma_1$ or $\gamma_2$ necessarily leads to better expectations because it simultaneously brings a large computation and communication burden.

For the relationship between the setting of the parameter $\gamma_1$ (or $\gamma_2$) and the regret bound, $\textbf{N}_{LO}$, communication number, Table \ref{tradeoff table} shows in detail some examples.
Taking the setting $\gamma_2=1$ in Table \ref{tradeoff table} as an analysis example, the regret bound $\mathcal{O} (1+H_T)$ is achieved regardless of the order of $H_T$, which is undoubtedly a great bound from the perspective of algorithm performance. However, this parameter setting incurs the computation and communication burden in $\mathcal{O} (T^2)$, which deviates from our expectations. Further, if the order of $H_T$ over time $T$ is large (e.g., $\mathcal{O}(\sqrt{T})$ in Table \ref{tradeoff table}), this setting $\gamma_2=1$ may be unreasonable for improving performance and saving resources.

Moreover, when communication and computing resources are limited, a small parameter $\gamma_1$ or $\gamma_2$ is more suitable for Algorithm \ref{algorithm 2}.  Therefore, it is better that the parameters $\varepsilon_1, \gamma_1$ (or $\varepsilon_2, \gamma_2$)  satisfy the trade-off between convergence accuracy, computational cost, and communication cost and in practical applications.

  \begin{table}[ht] 
\renewcommand\arraystretch{1.3}
\begin{center}
  \caption{The relationship between  the parameter $\gamma_1$ ( or $\gamma_2$), the regret bound, $\textbf{N}_{LO}$, and communication number.}
 \resizebox{\linewidth}{!}{ 
\label{tradeoff table}
 \begin{tabular}{ccccc}
\hline
$H_T (\hat{H}_T)$ & $\gamma_1$ (or $\gamma_2$)  & $\textbf{N}_{LO}$&  \makecell[c]{Communication \\ number } &    $\textbf{Regret}_d^j(T)$\\ \hline 
 \multirow{4}{*}{Unknown}  & $0.3 $ &  $\mathcal{O} (T^{1.3})$  &  $\mathcal{O} (T^{1.3})$ & $\mathcal{O}(T^{0.7}+H_T)$\\
                           & $0.5 $  &  $\mathcal{O} (T^{1.5})$  &  $\mathcal{O} (T^{1.5})$ &$\mathcal{O}(\sqrt{T}+H_T)$ \\
                           & $0.7 $  &  $\mathcal{O} (T^{1.7})$  &  $\mathcal{O} (T^{1.7})$ &$\mathcal{O}(T^{0.3}+H_T)$\\
                           & $\gamma_2=1 $  & $\mathcal{O} (T^{2})$  & $\mathcal{O} (T^{2})$ &$\mathcal{O}(1+H_T)$\\ \hline
    \multirow{4}{*}{$\mathcal{O}(\sqrt{T})$}  & $0.3 $ &  $\mathcal{O} (T^{1.3})$  &  $\mathcal{O} (T^{1.3})$ & $\mathcal{O}(T^{0.7})$\\
                           & $0.5 $  &  $\mathcal{O} (T^{1.5})$  &  $\mathcal{O} (T^{1.5})$ &$\mathcal{O}(\sqrt{T})$ \\
                           & $0.7 $  &  $\mathcal{O} (T^{1.7})$  &  $\mathcal{O} (T^{1.7})$ &$\mathcal{O}(1+H_T)$\\
                           &$\gamma_2=1 $  & $\mathcal{O} (T^{2})$  & $\mathcal{O} (T^{2})$ &$\mathcal{O}(1+H_T)$\\ \hline
\end{tabular}
}
\end{center}
\end{table}
\section{Simulation}
The distributed ridge regression problem shown in (\ref{Sim_eq}) is investigated in this section to validate the performance of the proposed algorithm.
\begin{flalign} \label{Sim_eq}
&\underset{\boldsymbol{x}}{\text{minimize}} \quad \quad  \sum_{t=1}^T \sum_{i=1}^n\left[\frac{1}{2}\left(\boldsymbol{a}_{i, t}^{\top} \bsx -b_{i, t}\right)^2+\lambda_1\|\bsx\|_2^2\right] \nonumber \\
&\text{subject to} \ \ \quad  \boldsymbol{x}\in \boldsymbol{X}
\end{flalign}
where $\boldsymbol{a}_{i, t}  \in \mathbb{R}^d$ is the feature vector and generated randomly and  uniformly in $[-5, 5]^d$, and $ b_{i, t} \in \mathbb{R}$ represents the label information.   The label $b_{i, t}$ satisfies
\begin{flalign} \label{Sim_eq2}
 b_{i, t}=\boldsymbol{a}_{i, t}^{\top} \bsx_0 + \frac{\zeta_{i,t}}{4 {t}}
\end{flalign}
 where $\zeta_{i,t}$ is generated randomly in $[0, 1]$. In this simulation, we consider the following two constraints: (i) the unit simplex constraint $\boldsymbol{X}: = \{ \bsx | \boldsymbol{1}^T \bsx =1 \}, [\bsx]_i\geq 0$; (ii) the norm ball constraint $\boldsymbol{X}: = \{ \bsx | \| \bsx\|_1 \leq 2 \}$.
 To intuitively verify the convergence performance of the developed algorithm, we define the global average dynamic regret $\frac{1}{n}\sum_{j=1}^n [\textbf{Regret}_d^j(T)/T]$, the upper envelope $\text{sup}_j\{ \textbf{Regret}_d^j(T)/T\}$ and the lower envelope $\text{inf}_j\{ \textbf{Regret}_d^j(T)/T\}$ of average dynamic regret, respectively. In the following cases, we set the parameters $n=20$ and $\lambda_1=5 \times 10^{-6}$.

\subsection{The Unit Simplex Constraint}
\begin{figure}[t]
\centering
\includegraphics[width=8cm]{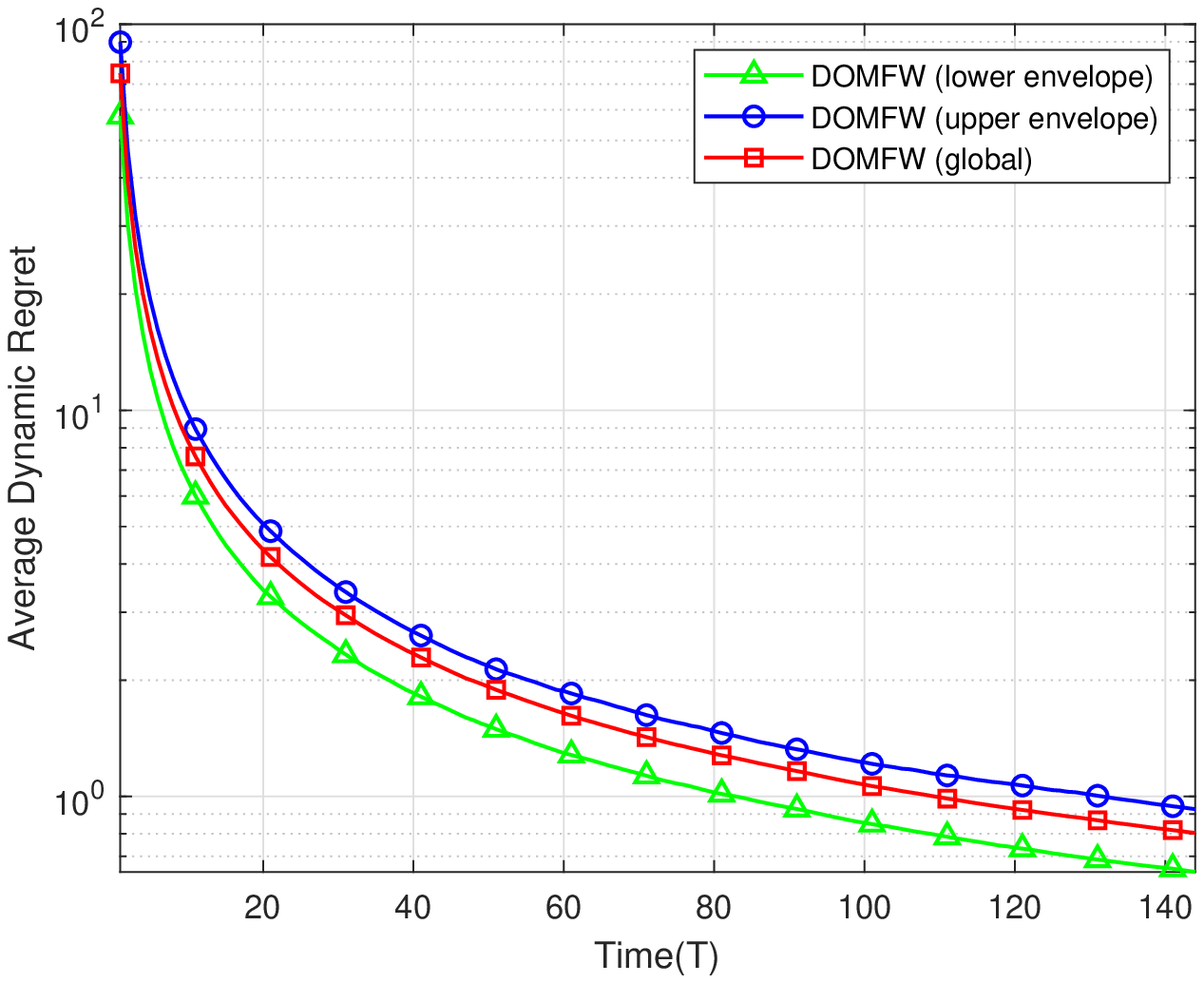}
 \captionsetup{font=small}
  \caption{Three average dynamic  regrets under the unit simplex constraint.}
  \label{MFW_1}
  \centering
\includegraphics[width=8cm]{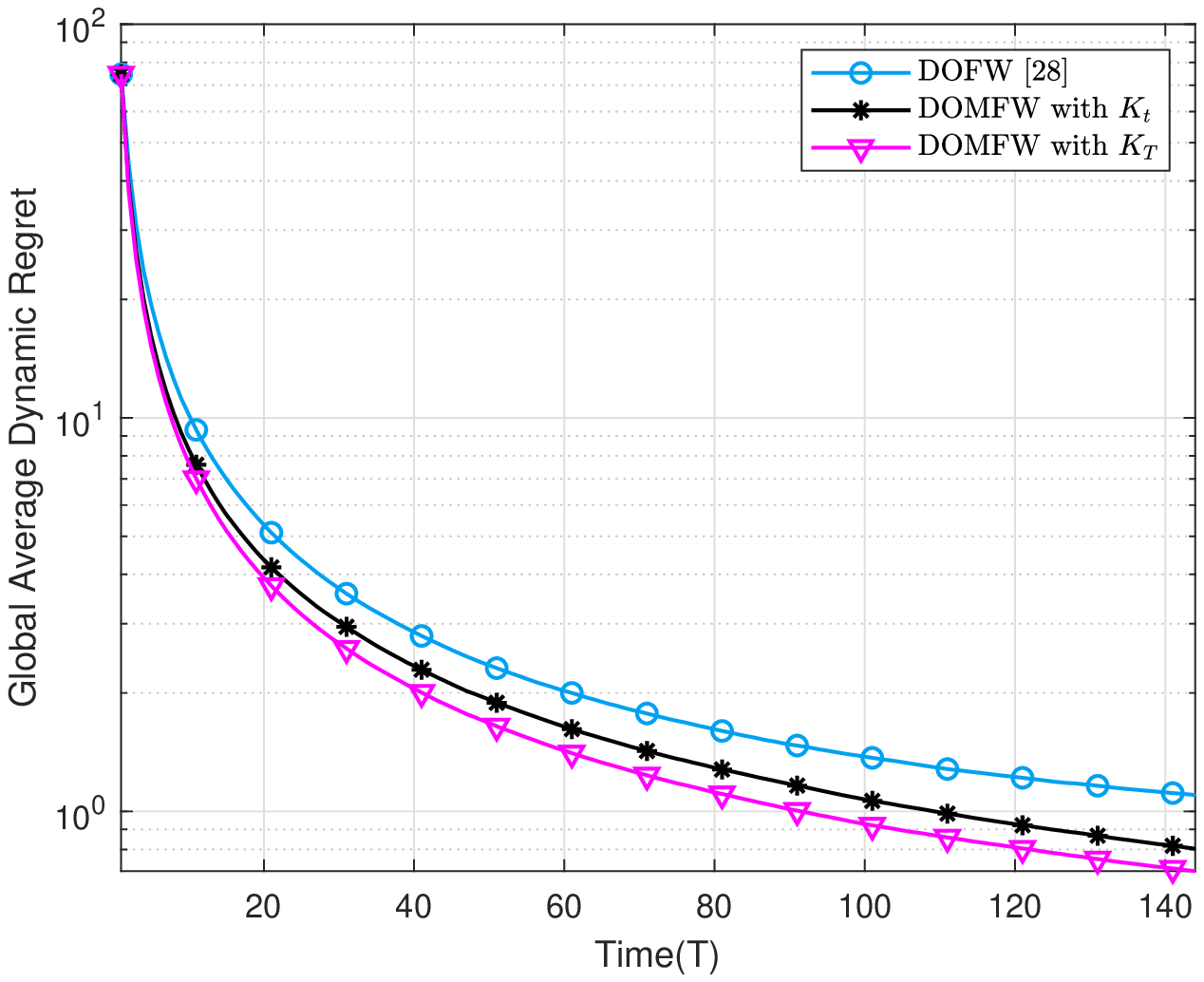}
 \captionsetup{font=small}
  \caption{ The comparison results of multiple iterations  under the unit simplex constraint.}
  \label{MFW_2}
    \centering
\includegraphics[width=8cm]{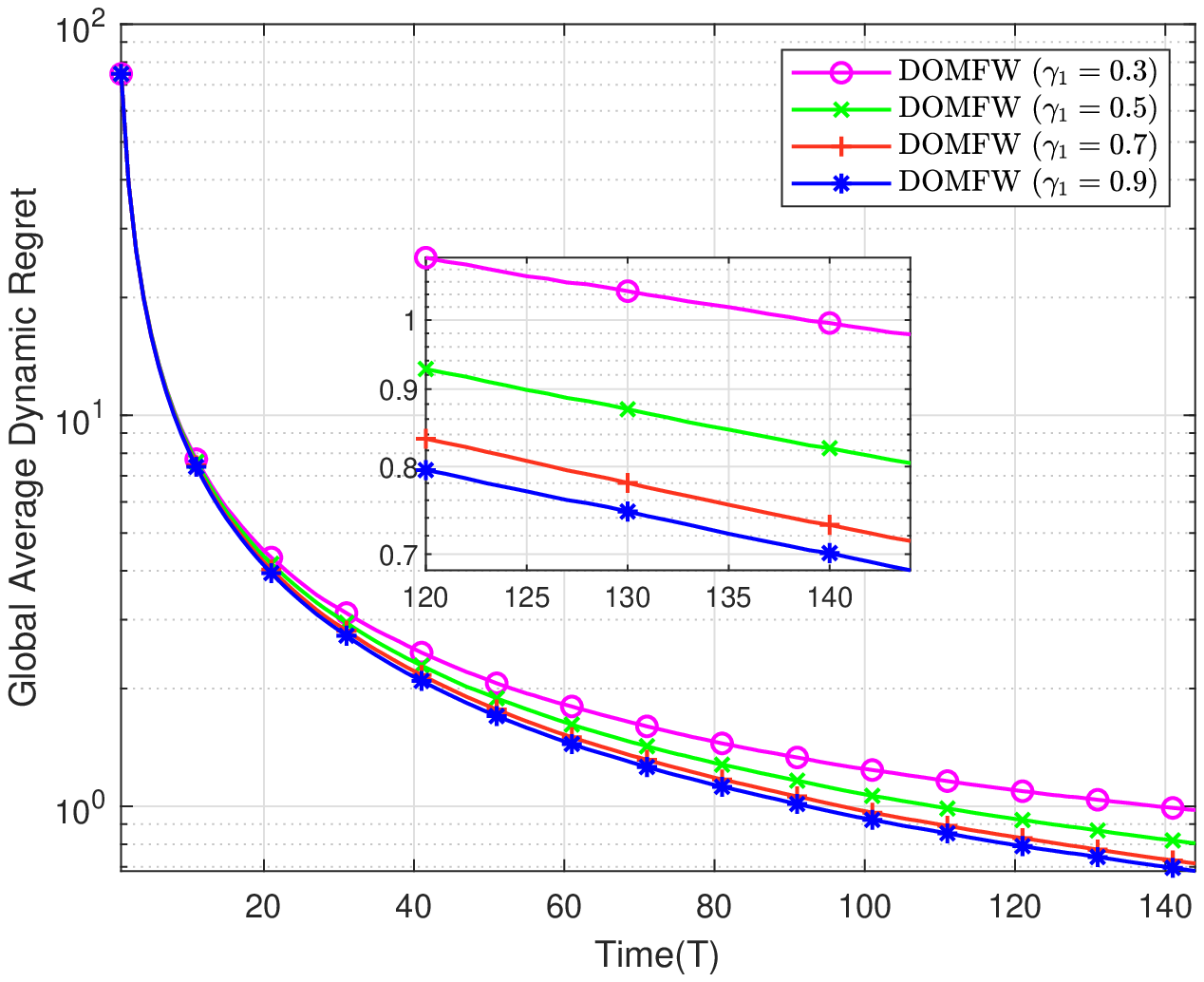}
 \captionsetup{font=small}
  \caption{ The comparison of global average dynamic regret for the unit simplex constraint under different inner number $K_t$.}
  \label{MFW_3}
\end{figure}

Under the unit simplex constraint, the convergence performance of Algorithm \ref{algorithm 2} is firstly analyzed. In Fig. \ref{MFW_1},  the simulation results show that the upper envelope, the lower envelope and the global value of $ \textbf{Regret}_d^j(T)/T$ are convergent for Algorithm \ref{algorithm 2} under the condition $d=8, K_t=\lceil 4 \sqrt{t} \rceil +1, \rho=4$, which is consistent with our theoretical results.
  Next, we compare the convergence performance of Algorithm \ref{algorithm 2} with the time-varying parameter $K_t$, with fixed parameter $K_T$ and distributed online Frank-Wolfe (DOFW) algorithm in \cite{zhang2023dynamic}, where the related parameters  are set as $ K_t=\lceil 4 \sqrt{t} \rceil +1, \rho=4$, $K_T=\lceil 4 \sqrt{T} \rceil +1, \rho=4$, and $\alpha=$$1/(4T^{0.4})$ \cite{zhang2023dynamic}, respectively. From Fig. \ref{MFW_2}, it is not hard to find  that the performance of Algorithm DOMFW is better than  that without the multiple iterations, which is consistent with the analysis in Remark \ref{corollary remark}(iii). Moreover, under this constraint, the algorithm with fixed parameter $K_T$ performs better than that with $K_t$.

 Now, we studies the effect of inner iteration number $K_t$ on the convergence performance of Algorithm \ref{algorithm 2} under the parameters  $\rho=3, \varepsilon_1=2$. The plots in Fig. \ref{MFW_3} show the global average dynamic regret for three different settings of $K_T$ in the same network, i.e., the cases when $\gamma_2$ is chosen as $0.3, 0.5, 0.7, 0.9$. The simulation results clearly reveal that the convergence performance of Algorithm \ref{algorithm 2} is getting better and better as $K_t$ increases, which corresponds to the results in Corollary \ref{corollary 3}. Meanwhile, when the parameter $\gamma_1$ changes from $0.7$ to $0.9$, the improvement of the convergence effect becomes slow.  From (\ref{Sim_eq2}), the estimation $\hat{H}_T$ can be approximated as $\mathcal{O}(\ln T)$. Thus, as the parameter $\gamma_1$ gets larger, especially for $\gamma_1\in [0.7,1)$, the effect of $H_T$ in the  upper bound of $\textbf{Regret}_d^j(T)$ is more emphasized and the performance improvement brought by $\gamma_1$ is weaker.

 \subsection{The Norm Ball Constraint}

\begin{figure}[t]
  \centering
\includegraphics[width=8cm]{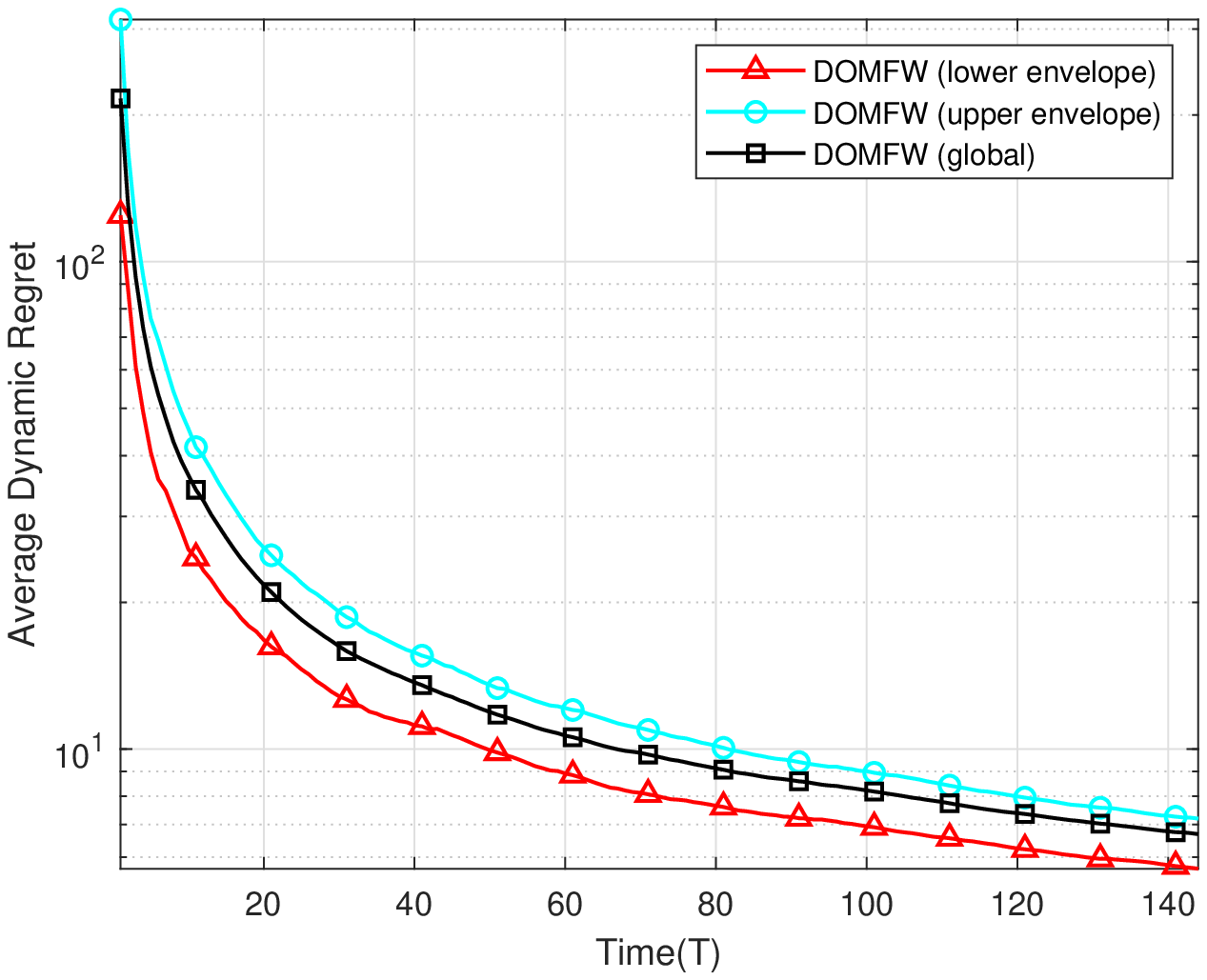}
 \captionsetup{font=small}
  \caption{Three average dynamic  regrets of Algorithm \ref{algorithm 2} under the norm ball constraint.}
  \label{MFW_5}
  \centering
\includegraphics[width=8cm]{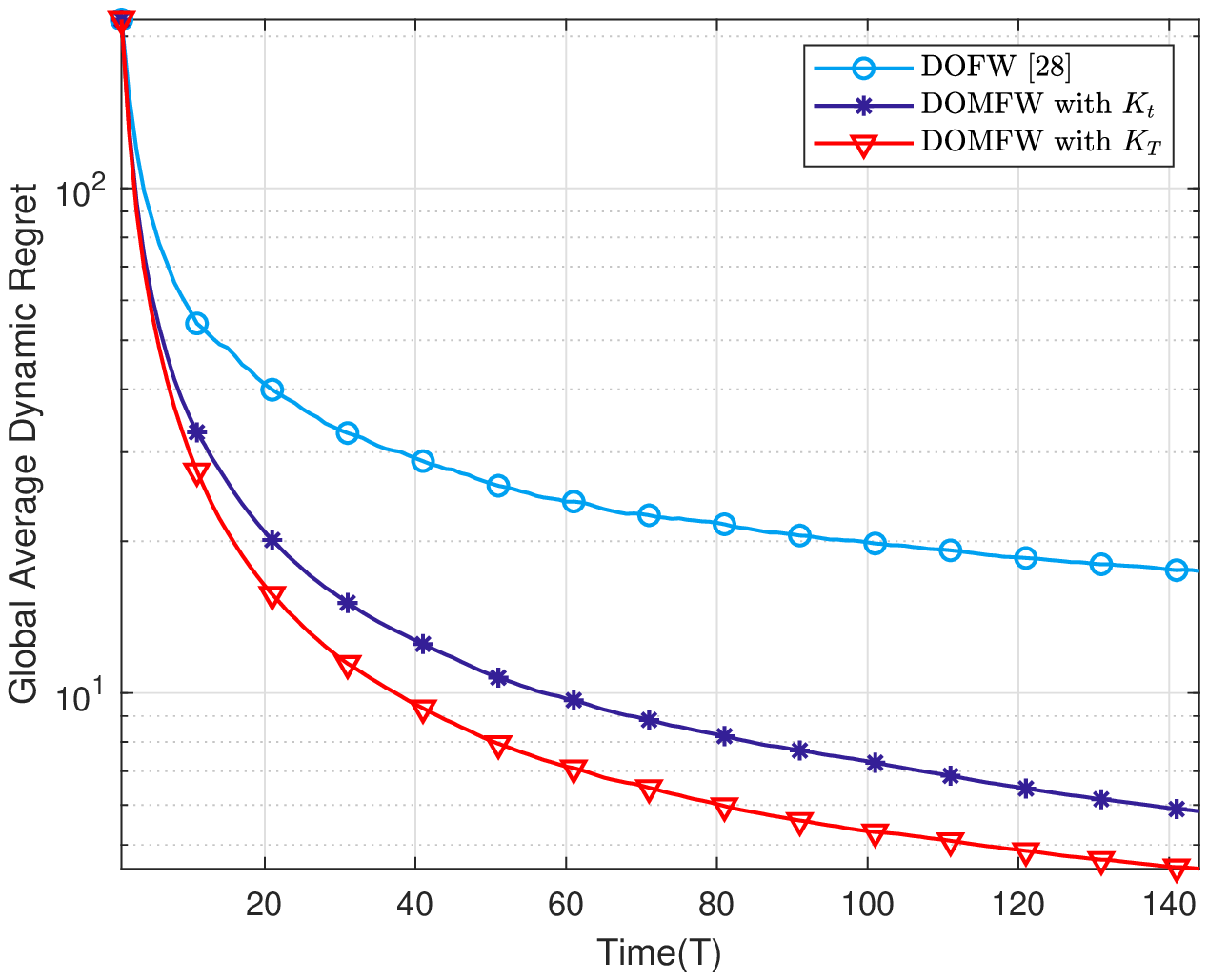}
 \captionsetup{font=small}
  \caption{ The comparison results of multiple iterations  under the norm ball constraint.}
  \label{MFW_4}
    \centering
\includegraphics[width=8cm]{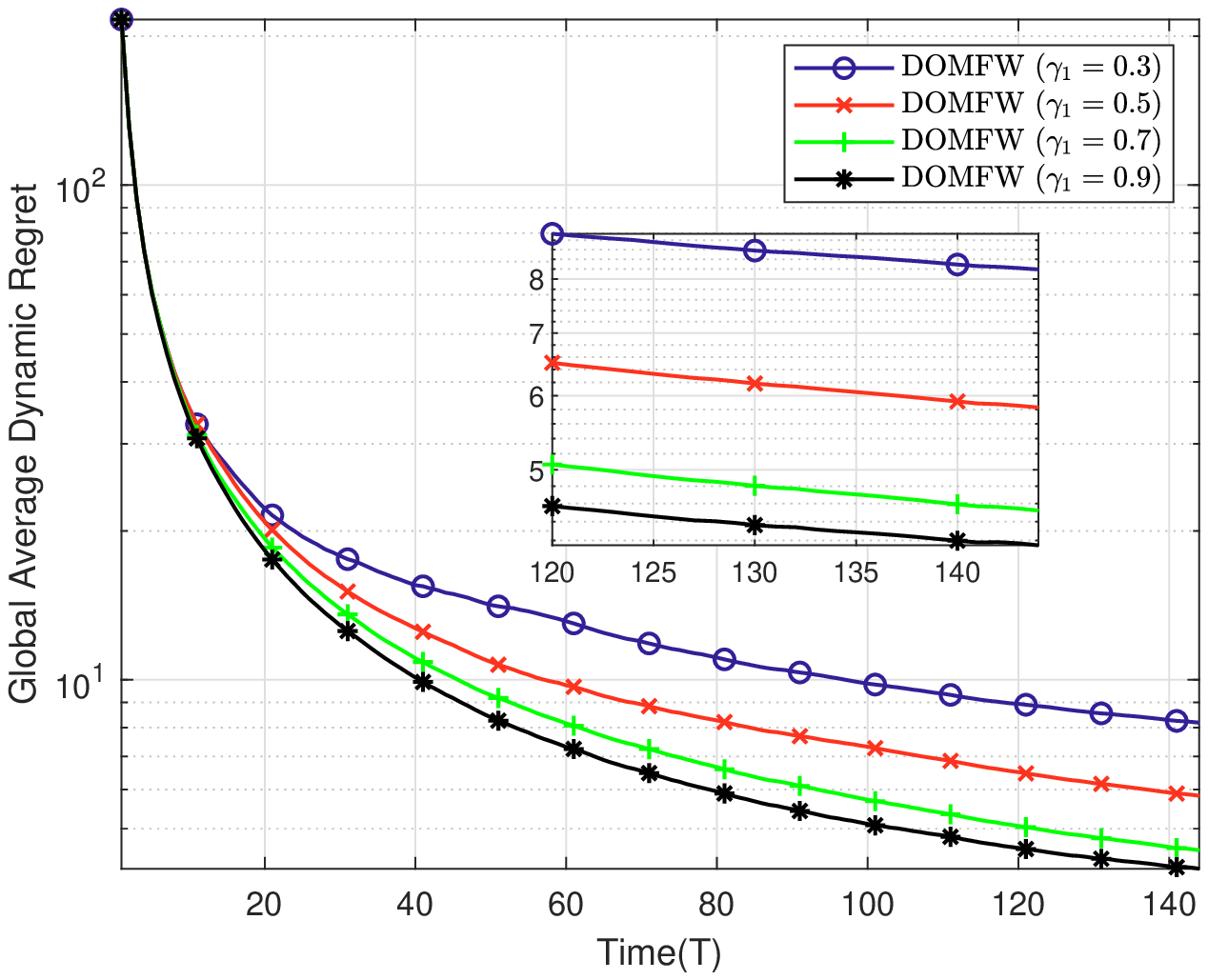}
 \captionsetup{font=small}
  \caption{  The comparison of global average dynamic regret for the norm ball constraint under different inner number $K_t$.}
  \label{MFW_6}
\end{figure}
 In this subsection, we study the convergence performance of Algorithm \ref{algorithm 2} under  the norm ball constraint with the dimension $d=16$.
  Set $ \varepsilon_1=\varepsilon_2=0.3,$ $ \gamma_1=\gamma_2=0.5, \rho=3.5$. Similarly, by observing three regret curves in Fig. \ref{MFW_4}, we obtain that Algorithm \ref{algorithm 2} is convergent for this constraint.
  Next, to verify the effect of the multiple iteration method, we compared  DOFW algorithm in \cite{zhang2023dynamic} and Algorithm \ref{algorithm 2} under two inner iteration parameters.   From Fig. \ref{MFW_4}, Algorithm DOMFW has better convergence performance than that without multiple iterations \cite{zhang2023dynamic}, which reflects the significant advantages of the multiple iteration method. For two parameter settings $K_t$ and $K_T$, the performance generated from the latter is better, while the former does not depend on prior knowledge of $T$.

  Finally, taking $K_t$ as an example, we explore the convergence performance of Algorithm \ref{algorithm 2} under three different settings of $\gamma_1$, i.e., $0.3, 0.5, 0.7, 0.9$ and the settings $\rho=3.5, \varepsilon_1=0.3$.
  Fig. \ref{MFW_6} clearly shows that as the order setting of $K_t$ over time $t$ increases, the global average dynamic regret of Algorithm \ref{algorithm 2} is getting better.
  Similar to Fig. \ref{MFW_3},  when $\gamma_1$ changes from $0.7$ to $0.9$, the improvement of the convergence effect is weak while the increase in the computational cost and communication cost is great. Therefore, as described in the previous section, a proper setting of parameters $\varepsilon_1 (\text{or} \ \varepsilon_2)$ and $\gamma_1 (\text{or} \ \gamma_2)$ is very important to trade off the relationship between obtaining high-quality decision and saving resources in practical applications.
\section{Conclusions}
 For solving the  distributed online optimization problem, we have developed the distributed online multiple Frank-Wolfe  algorithm   over a time-varying multi-agent network. Based on the projection-free operation, the proposed algorithm can significantly save computational cost (at one step), especially for the optimization scenario with a high-dimensional and  structural constraint set. Further, we have confirmed that the multiple iteration technique can enhance the dynamic regret bound of the FW algorithm in distributed scenarios.
The regret bound $\mathcal{O}(T^{1-\gamma}+H_T),0<\gamma<1$ with the linear oracle number $\mathcal{O}(T^{1+\gamma})$ has been established, which is tighter than that without inner iteration loop and does not require a step size dependent on the prior knowledge of $H_T,T$. In particular, we have confirmed that when the order or estimated order  of $H_T$ over time $T$  is available, the optimal bound $\mathcal{O} (1+H_T)$  can be obtained. Moreover, a trade-off between dynamic regret bound, computational cost, and communication cost has been revealed. Finally, we have verified the theoretical results by the simulation of the distributed online ridge regression problems with two constraint sets.

\begin{appendices}
\section{Proof of Lemma \ref{consistency2}}

According to Algorithm \ref{algorithm 2}, by denoting $A_t^{K_t}$ as the product of ${K_t}$ multiplications of $A_t$, we get
\begin{flalign} \label{proof_consistency2 1}
&{\bsx}_{i,t+1} \nonumber \\
&= \hat{\boldsymbol{x}}_{i,t}^{K_t} +\alpha_t (\boldsymbol{v}_{i,t}^{K_t}- \hat{\boldsymbol{x}}_{i,t}^{K_t})\nonumber \\
 &= \sum\limits_{j=1}^n [A_t]_{ij}  {\bsx}_{j,t}^{K_t} +\alpha_t  \left(\boldsymbol{v}_{i,t}^{K_t}-\hat{\bsx}_{i,t}^{K_t} \right)\nonumber \\
%
%
&= \sum_{j=1}^n [A_t^{K_t}]_{ij}   {\bsx}_{j,t} +\alpha_t \sum\limits_{l=1}^{{K_t}-1} \sum\limits_{j=1}^n [A_t^{{K_t}-l}]_{ij} \left(\boldsymbol{v}_{j,t}^l-\hat{\bsx}_{j,t}^l \right)\nonumber \\
&\quad + \alpha_t \left(\boldsymbol{v}_{i,t}^{K_t}-\hat{\bsx}_{i,t}^{K_t} \right)\nonumber \\
&= \sum_{j=1}^n [\Phi^K(t,1)]_{ij}   {\bsx}_{j,1} + \alpha_t \left(\boldsymbol{v}_{i,t}^{K_t}-\hat{\bsx}_{i,t}^{K_t} \right)          \nonumber \\
&\quad +\sum_{s=1}^{t}\sum\limits_{l=1}^{{K_{s}}-1} \sum\limits_{j=1}^n \alpha_{s}[\Phi^K(t,s+1) A_{s}^{{K_{s}}-l}  ]_{ij} \left(\boldsymbol{v}_{j,s}^l-\hat{\bsx}_{j,s}^l \right) \nonumber \\
&\quad + \sum_{s=1}^{t-1} \sum_{j=1}^n  \alpha_s[\Phi^K(t,s+1)   ]_{ij}\left(\boldsymbol{v}_{j,s}^{K_{s}}-\hat{\bsx}_{j,s}^{K_{s}} \right).
%
%
\end{flalign}

According to Algorithm \ref{algorithm 2}, we have that
\begin{flalign} \label{proof_consistency2 2-a1}
& \quad {\boldsymbol{x}}_{avg,t}^{K_t+1}\nonumber \\
&= \frac{1}{n} \sumn \bsx_{i,t}^{{K_t}+1} \nonumber \\
 &=\frac{1}{n} \sumn \left[\sum\limits_{j=1}^n [A_t]_{ij}  {\bsx}_{j,t}^{K_t} +\alpha_t  \left(\boldsymbol{v}_{i,t}^{K_t}-\sum\limits_{j=1}^n [A_t]_{ij}  {\bsx}_{j,t}^{K_t}  \right)\right]\nonumber \\
&={\boldsymbol{x}}_{avg,t}^{K_t} +\alpha_t ({\boldsymbol{v}}_{avg,t}^{K_t}-{\boldsymbol{x}}_{avg,t}^{K_t}).
\end{flalign}
By iteration with respect to $K_t$, we further obtain that
\begin{flalign} \label{proof_consistency2 2-a2}
{\boldsymbol{x}}_{avg,t}^{K_t+1}
%
&={\boldsymbol{x}}_{avg,t} +\alpha_t \sum\limits_{l=1}^{{K_t}} ({\boldsymbol{v}}_{avg,t}^l-{\boldsymbol{x}}_{avg,t}^l).
\end{flalign}
Note that $
{\boldsymbol{x}}_{avg,t+1} ={\boldsymbol{x}}_{avg,t}^{K_t+1}.
$
 Therefore, we have
\begin{flalign} \label{proof_consistency2 2-a4}
&\quad {\boldsymbol{x}}_{avg,t+1}\nonumber \\
&={\boldsymbol{x}}_{avg,t} +\alpha_t \sum\limits_{l=1}^{{K_t}} ({\boldsymbol{v}}_{avg,t}^l-{\boldsymbol{x}}_{avg,t}^l)\nonumber \\
&={\boldsymbol{x}}_{avg,1} + \sum_{s=1}^t \sum\limits_{l=1}^{{K_s}} \alpha_{s} ({\boldsymbol{v}}_{avg,s}^l-{\boldsymbol{x}}_{avg,s}^l) \nonumber \\
%
&=\frac{1}{n}\sumjn {\bsx}_{j,1} +\frac{1}{n} \sum_{s=1}^t\sum\limits_{l=1}^{{K_{s}}-1}\sumjn  \alpha_{s}({\boldsymbol{v}}_{j,s}^l-\hat{\boldsymbol{x}}_{j,s}^l)\nonumber \\
&\quad +\frac{1}{n} \sum_{s=1}^{t-1}\sumjn \alpha_{s}({\boldsymbol{v}}_{j,s}^{K_{s}}-\hat{\boldsymbol{x}}_{j,s}^{K_{s}})+ \frac{\alpha_t}{n} \sumjn ({\boldsymbol{v}}_{j,t}^{K_t}-\hat{\boldsymbol{x}}_{j,t}^{K_t}).
\end{flalign}

By combining (\ref{proof_consistency2 1}) and  (\ref{proof_consistency2 2-a4}), it follows that
\begin{flalign} \label{proof_consistency2 3-a1}
&\quad \| {\bsx}_{i,t+1}- \bsx_{avg,t+1}\| \nonumber \\
 &\leq   \sum_{j=1}^n \left|[\Phi^K(t,1)]_{ij}- \frac{1}{n}  \right|\left\| {\bsx}_{j,1}\right\|  + \alpha_t \left\|\boldsymbol{v}_{i,t}^{K_t}-\hat{\bsx}_{i,t}^{K_t} \right\| + \nonumber \\
 &\sum_{s=1}^{t}\sum\limits_{l=1}^{{K_{s}}-1} \sum\limits_{j=1}^n \alpha_{s}\left|[\Phi^K(t,s+1) A_{s}^{{K_{s}}-l}  ]_{ij}  -\frac{1}{n}\right| \left\|\boldsymbol{v}_{j,s}^l-\hat{\bsx}_{j,s}^l \right\| \nonumber \\
&\quad  + \sum_{s=1}^{t-1} \sum_{j=1}^n  \alpha_s\left|[\Phi^K(t,s+1)   ]_{ij} -\frac{1}{n}\right|\left\|\boldsymbol{v}_{j,s}^{K_{s}}-\hat{\bsx}_{j,s}^{K_{s}} \right\|  \nonumber \\
&\quad +\frac{\alpha_t}{n} \sumjn \left\|{\boldsymbol{v}}_{j,t}^{K_t}-\hat{\boldsymbol{x}}_{j,t}^{K_t}\right\|.
\end{flalign}

Note  that, for any appropriate numbers $j,s,l$, the following inequality holds:
 \begin{flalign} \label{proof_consistency2 3-a2}
\left\|{\boldsymbol{v}}_{j,s}^{l}-\hat{\boldsymbol{x}}_{j,s}^{l}\right\|\leq M.
\end{flalign}
Then, through recalling (\ref{phi-K-b2}),  (\ref{phi-K-b3}) and using (\ref{proof_consistency2 3-a2}), (\ref{proof_consistency2 3-a1}) follows  that
\begin{flalign} \label{proof_consistency2 3-a3}
 & \quad \| {\bsx}_{i,t+1}- \bsx_{avg,t+1}\| \nonumber \\
&\leq  \Gamma_1 \sigma_1^{\sum_{p=1}^tK_p-1} \sum_{j=1}^n \left\| {\bsx}_{j,1}\right\|+nM \Gamma_1 \sum_{s=1}^{t-1}   \alpha_s \sigma_1^{\sum_{p=s+1}^tK_p-1}  \nonumber \\
&\quad   + nM\Gamma_1\sum_{s=1}^{t}\sum\limits_{l=1}^{{K_{s}}-1}  \alpha_{s} \sigma_1^{\sum_{p=s}^tK_p-l-1} + 2\alpha_t M \nonumber\\
&\leq  \Gamma_1 \sigma_1^{tK_1-1} \sum_{j=1}^n \left\| {\bsx}_{j,1}\right\|+nM \Gamma_1 \sum_{s=1}^{t}   \alpha_s \sigma_1^{(t-s)K_1-1}  \nonumber \\
&\quad + nM\Gamma_1\sum_{s=1}^{t}\sum\limits_{l=1}^{{K_{s}}-1}  \alpha_{s} \sigma_1^{(t-s)K_1+K_t-l-1}  + 2\alpha_t M
\end{flalign}
where in the last inequality we use the non-decreasing property of the sequence $\{K_t\}$.
Based on (\ref{proof_consistency2 3-a3}), we can further obtain that
\begin{flalign} \label{proof_consistency2 3-a3}
&\quad\sumT \sumn \| {\bsx}_{i,t}- \bsx_{avg,t}\| \nonumber \\
 &= \sum_{t=1}^{T-1} \sumn \| {\bsx}_{i,t+1}- \bsx_{avg,t+1}\|+ \sumn \| {\bsx}_{i,1}- \bsx_{avg,1}\| \nonumber \\
 &\leq n\Gamma_1   \sum_{t=1}^{T-1}  \sigma_1^{tK_1-1} \sum_{j=1}^n \left\| {\bsx}_{j,1}\right\|+ \sumn \| {\bsx}_{i,1}- \bsx_{avg,1}\| \nonumber \\
 &\quad + n^2M\Gamma_1 \sum_{t=1}^{T-1} \sum_{s=1}^{t}\sum\limits_{l=1}^{{K_{s}}-1}  \alpha_{s} \sigma_1^{(t-s)K_1+K_t-l-1} \nonumber \\
&\quad  +n^2M \Gamma_1  \sum_{t=1}^{T-1} \sum_{s=1}^{t}   \alpha_s \sigma_1^{(t-s)K_1-1}  + 2n M\sum_{t=1}^{T-1}\alpha_t.
\end{flalign}

Now, we are going to handle the summations involved in (\ref{proof_consistency2 3-a3}). Recalling that $0<\sigma_1<1$ and $K_1\geq 2$, we have
\begin{flalign}\label{proof_consistency2 6_a1}
  \sum_{t=1}^{T-1}  \sigma_1^{tK_1-1}
=  \frac{1}{ \sigma_1}\sum_{t=1}^{T-1}  \left(\sigma_1^{K_1}\right)^t \leq\frac{\sigma_1^{K_1}}{\sigma_1 \left(1-\sigma_1^{K_1}\right)}
\leq  \frac{1}{ 1-\sigma_1^{K_1} }.
\end{flalign}
Note that the sequence $\{K_t\}$ is non-decreasing. Then, it is not difficult to obtain that
\begin{flalign}&\quad \sum_{t=1}^{T-1} \sum_{s=1}^{t}\sum\limits_{l=1}^{{K_{s}}-1}  \alpha_{s} \sigma_1^{(t-s)K_1+K_t-l-1}  \nonumber \\
& \leq\sum_{t=1}^{T-1} \sum_{s=1}^{t} \alpha_{s} \sigma_1^{(t-s)K_1}\sum\limits_{l=1}^{{K_{t}}-1}\sigma_1^{K_t-l-1}  \nonumber \\
&\leq \frac{1}{1-{\sigma_1}}\sum_{t=1}^{T-1}  \sum_{s=1}^{t}\alpha_{s} {\sigma_1}^{(t-s){K_{1}}}  \nonumber \\
&\leq\frac{1}{1-{\sigma_1}}\left(\sum_{t=1}^{T}  \alpha_{t}\right)\left(\sum_{s=1}^{T} {\sigma_1}^{(s-1){K_{1}}}\right)  \nonumber \\
&\leq \frac{1}{\left(1-{\sigma_1}\right)\left(1-{\sigma_1}^{K_1}\right) } \sumT  \alpha_{t}.  \label{proof_consistency2 6_a2}    
\end{flalign}
Similarly, it can be verified that
\begin{flalign}
\sum_{t=1}^{T-1} \sum_{s=1}^{t}   \alpha_s \sigma_1^{(t-s)K_1-1}
& \leq\frac{1}{{\sigma_1} \left(1-{\sigma_1}^{K_1}\right) } \sumT  \alpha_{t}.   \label{proof_consistency2 6_a3}
\end{flalign}
Combining (\ref{proof_consistency2 3-a3})-(\ref{proof_consistency2 6_a3}) yields the condition in (\ref{consistency2-con1}). The proof is complete. \hfill$\square$


\section{Proof of Lemma \ref{grad diffience2}}
 Under Assumption \ref{assump: lips grad}, we have
\begin{flalign}\label{proof_grad diff2 1}
 &\quad \left\|  \widehat{\nabla} f_{i,t}^k - \frac{1}{n}\nabla F_t ({\bsx}_{avg,t}^k)\right\|\nonumber \\
 &\leq  \left\|  \widehat{\nabla} f_{i,t}^k - \frac{1}{n}\sumn \nabla f_{i,t} (\hat{\bsx}_{i,t}^k)\right\|\nonumber \\
  &\quad + \left\|   \frac{1}{n}\sumn \nabla f_{i,t} (\hat{\bsx}_{i,t}^k) - \frac{1}{n}\nabla F_t ({\bsx}_{avg,t}^k)\right\| \nonumber \\
  &\leq  \left\|  \widehat{\nabla} f_{i,t}^k - \frac{1}{n}\sumn \nabla f_{i,t} (\hat{\bsx}_{i,t}^k)\right\| + \frac{G_X}{n}\sumn \left\|  \hat{\bsx}_{i,t}^k -  {\bsx}_{avg,t}^k\right\| \nonumber \\
   &\leq  \left\|  \widehat{\nabla} f_{i,t}^k - \frac{1}{n}\sumn \nabla f_{i,t} (\hat{\bsx}_{i,t}^k)\right\| + \frac{G_X}{n}\sumn \left\|  {\bsx}_{i,t}^k -  {\bsx}_{avg,t}^k\right\|
\end{flalign}
where the last inequality is obtained by using the double stochasticity of $A_t$. Then, after summing (\ref{proof_grad diff2 1}) with respect to $i, k$ and $t$, we have
\begin{flalign}\label{proof_grad diff2 2}
  &\quad \sumT \sum_{k=1}^{K_t} \sumn \alpha_t\left\|  \widehat{\nabla} f_{i,t}^k - \frac{1}{n}\nabla F_t ({\bsx}_{avg,t}^k)\right\| \nonumber \\
   &\leq   \sumT \sum_{k=1}^{K_t} \sumn \alpha_t \left\|  \widehat{\nabla} f_{i,t}^k - \frac{1}{n}\sumn \nabla f_{i,t} (\hat{\bsx}_{i,t}^k)\right\|  \nonumber \\
   &\quad +  G_X\sumT \sum_{k=1}^{K_t} \sumn \alpha_t\left\|  {\bsx}_{i,t}^k -  {\bsx}_{avg,t}^k\right\|.
\end{flalign}

The first term on the right hand side  is first analysed.
\begin{flalign} \label{proof_grad diff2 3}
& \quad \widehat{\nabla} f_{i,t}^k \nonumber \\
&=\sum\limits_{j=1}^n [A_t]_{ij}  \overline{\nabla} f_{j,t}^k\nonumber \\
%
%
&= \sum\limits_{j=1}^n [A_t]_{ij}\widehat{\nabla} f_{j,t}^{k-1} + \sum\limits_{j=1}^n [A_t]_{ij}  \delta_{j,t}^{k-1} \nonumber \\
%
%
&= \sum\limits_{j=1}^n [A_t^k]_{ij}{\nabla} f_{j,t} (\hat{\bsx}_{j,t}^1) + \sum\limits_{l=2}^k \sum\limits_{j=1}^n [A_t^{k-l+1}]_{ij} \delta_{j,t}^{l-1}.
\end{flalign}

Based on Algorithm \ref{algorithm 2}, it is easily obtained that $ \sum_{i=1}^n \overline{\nabla} f_{i,1}^1 =   \sum_{i=1}^n  \nabla f_{i,1} (\hat{\boldsymbol{x}}_{i,1}^1)$. Now we assume that  $\sum_{i=1}^n \overline{\nabla} f_{i,t}^{k-1} = \sum_{i=1}^n \nabla f_{i,t} (\hat{\boldsymbol{x}}_{i,t}^{k-1})$ holds. Then, we intend to prove the equality also hold at time $t+1$. Actually,
\begin{flalign} \label{proof_grad diff2 5}
&\quad \sumn \overline{\nabla} f_{i,t}^k \nonumber \\
& = \sumn \widehat{\nabla} f_{i,t}^{k-1} + \sumn  \nabla f_{i,t}  (\hat{\boldsymbol{x}}_{i,t}^k) - \sumn  \nabla f_{i,t}  (\hat{\boldsymbol{x}}_{i,t}^{k-1}) \nonumber \\
&= \sumn \sum_{j=1}^n [A_t]_{ij} \overline{\nabla} f_{j,t}^{k-1} + \sumn  \nabla f_{i,t}  (\hat{\boldsymbol{x}}_{i,t}^k) - \sum_{i=1}^n \overline{\nabla} f_{i,t}^{k-1}  \nonumber \\
&=\sumn  \nabla f_{i,t}  (\hat{\boldsymbol{x}}_{i,t}^k)
\end{flalign}
where the third equality combines the double stochasticity of $A_{t}$.

This further gives that
\begin{flalign} \label{proof_grad diff2 4}
&\quad \frac{1}{n} \sumn   {\nabla} f_{i,t} (\hat{\boldsymbol{x}}_{i,t}^k)  \nonumber \\
%
&= \frac{1}{n} \sumn \sum\limits_{j=1}^n [A_{t}]_{ij}  \overline{\nabla} f_{j,t}^{k-1} +  \frac{1}{n} \sumn  \delta_{i,t}^{k-1}  \nonumber \\
%
&= \frac{1}{n}\sumn  \nabla f_{i,t}  (\hat{\boldsymbol{x}}_{i,t}^{k-1})+  \frac{1}{n} \sumn  \delta_{i,t}^{k-1}    \nonumber \\
&= \frac{1}{n} \sumn   {\nabla} f_{i,t} (\hat{\boldsymbol{x}}_{i,t}^1) +  \frac{1}{n}\sum\limits_{l=2}^k \sumn   \delta_{i,t}^{l-1}.
\end{flalign}

Based on the above analysis, we obtain that, for any $k\geq2$,
\begin{flalign} \label{proof_grad diff2 6}
 &\quad \left\|  \widehat{\nabla} f_{i,t}^k - \frac{1}{n} \sumn   {\nabla} f_{i,t} (\hat{\boldsymbol{x}}_{i,t}^k)\right\|\nonumber \\
 &\leq\sum\limits_{j=1}^n \left | [A_t^k]_{ij} -  \frac{1}{n} \right|  \left\|   {\nabla} f_{j,t} (\hat{\boldsymbol{x}}_{j,t}^1) \right\| \nonumber \\
&\quad +\sum\limits_{l=2}^k \sum\limits_{j=1}^n \left | [A_t^{k-l+1}]_{ij} -  \frac{1}{n} \right| \left\|  \delta_{j,t}^{l-1} \right\|\nonumber \\
  &\leq \Gamma_1 {\sigma_1}^{k-1}\sum\limits_{j=1}^n \left\|   {\nabla} f_{j,t} (\hat{\boldsymbol{x}}_{j,t}^1) \right\| + \Gamma_1 \sum\limits_{l=2}^k \sum\limits_{j=1}^n{\sigma_1}^{k-l} \left\|  \delta_{j,t}^{l-1} \right\|.
\end{flalign}
This implies that
\begin{flalign} \label{proof_grad diff2 7}
 &\quad \sum_{k=1}^{K_t} \sumn \left\|  \widehat{\nabla} f_{i,t}^k - \frac{1}{n} \sumn   {\nabla} f_{i,t} (\hat{\boldsymbol{x}}_{i,t}^k)\right\|\nonumber \\
 &=\sumn \left\|  \widehat{\nabla} f_{i,t}^1 - \frac{1}{n} \sumn   {\nabla} f_{i,t} (\hat{\boldsymbol{x}}_{i,t}^1)\right\|  \nonumber \\
 &\quad +\sum_{k=2}^{K_t} \sumn \left\|  \widehat{\nabla} f_{i,t}^k - \frac{1}{n} \sumn   {\nabla} f_{i,t} (\hat{\boldsymbol{x}}_{i,t}^k)\right\| \nonumber \\
 &\leq\sumn \left\|  \widehat{\nabla} f_{i,t}^1 - \frac{1}{n} \sumn   {\nabla} f_{i,t} (\hat{\boldsymbol{x}}_{i,t}^1)\right\|   +n\Gamma_1 \sum\limits_{k=2}^{K_t} \sum\limits_{j=1}^n {\sigma_1}^{k-1} \cdot \nonumber \\
 &\quad \left\|   {\nabla} f_{j,t} (\hat{\boldsymbol{x}}_{j,t}^1) \right\| +n\Gamma_1\sum\limits_{k=2}^{K_t} \sum\limits_{l=2}^k \sum\limits_{j=1}^n {\sigma_1}^{k-l} \left\|  \delta_{j,t}^{l-1} \right\|.
\end{flalign}
Now we are going to estimate the bounds of the terms in  (\ref{proof_grad diff2 7}). Firstly, we have
\begin{flalign}
&\quad \sumn \left\|  \widehat{\nabla} f_{i,t}^1 - \frac{1}{n} \sumn   {\nabla} f_{i,t} (\hat{\boldsymbol{x}}_{i,t}^1)\right\| \nonumber \\
 &= \sumn \left\| \sum\limits_{j=1}^n [A_t]_{ij} {\nabla} f_{j,t} (\hat{\boldsymbol{x}}_{j,t}^1) - \frac{1}{n} \sum\limits_{j=1}^n \nabla f_{j,t} (\hat{\bsx}_{j,t}^1)\right\| \nonumber \\
&\leq \sumn \sum\limits_{j=1}^n \left | [A_t]_{ij} -  \frac{1}{n} \right|  \left\|  {\nabla} f_{j,t} (\hat{\boldsymbol{x}}_{j,t}^1) \right\|\nonumber \\
&\leq n^2 \Gamma_1 L_X. \label{proof_grad diff2 8_1}
\end{flalign}
Next, it is easily obtained that
\begin{flalign}
n\Gamma_1\sum\limits_{k=2}^{K_t} \sum\limits_{j=1}^n {\sigma_1}^{k-1}  \left\|   {\nabla} f_{j,t} (\hat{\boldsymbol{x}}_{j,t}^1) \right\|
 &\leq \frac{{\sigma_1} n^2 \Gamma_1 L_X}{1-{\sigma_1}}. \label{proof_grad diff2 8_2}\end{flalign}
Moreover, it can be verified that
\begin{flalign}
 &\quad n\Gamma_1\sum\limits_{k=2}^{K_t} \sum\limits_{l=2}^k \sum\limits_{j=1}^n {\sigma_1}^{k-l} \left\|  \delta_{j,t}^{l-1} \right\|\nonumber \\
  &\leq
n\Gamma_1 \left(\sum\limits_{k=2}^{{K_t} }
 {\sigma_1}^{k-2}\right)\left(\sum\limits_{k=2}^{{K_t} } \sumn \left\|  \delta_{i,t}^{k-1}  \right\|\right)\nonumber \\
  %
  &\leq \frac{ n \Gamma_1}{1-{\sigma_1}} \sum\limits_{k=1}^{{K_t}-1} \sumn \left\|  \nabla f_{i,t}  (\hat{\boldsymbol{x}}_{i,t}^{k+1}) - \nabla f_{i,t}  (\hat{\boldsymbol{x}}_{i,t}^{k}) \right\| \nonumber \\
  &\leq \frac{ n \Gamma_1 G_X}{1-{\sigma_1}} \sum\limits_{k=1}^{{K_t}-1} \sumn \left\|  \hat{\boldsymbol{x}}_{i,t}^{k+1} - \hat{\boldsymbol{x}}_{i,t}^{k} \right\|. \label{proof_grad diff2 8_3}
\end{flalign}

Note that
\begin{flalign} \label{proof_grad diff2 9}
 &\quad \sumn \left\|  \hat{\boldsymbol{x}}_{i,t}^{k+1} - \hat{\boldsymbol{x}}_{i,t}^{k} \right\|\nonumber \\
 &\leq   \sumn \left\|  \hat{\boldsymbol{x}}_{i,t}^{k+1} - {\boldsymbol{x}}_{avg,t}^{k} \right\| +\sumn \left\|  {\boldsymbol{x}}_{avg,t}^{k} - \hat{\boldsymbol{x}}_{i,t}^{k} \right\|\nonumber \\
 &\leq   \sumn \left\|  {\boldsymbol{x}}_{i,t}^{k+1} - {\boldsymbol{x}}_{avg,t}^{k} \right\| +\sumn \left\|  {\boldsymbol{x}}_{avg,t}^{k} - {\boldsymbol{x}}_{i,t}^{k} \right\|\nonumber \\
 &\leq   \sumn \left\|  \hat{\boldsymbol{x}}_{i,t}^{k} +\alpha_t (\boldsymbol{v}_{i,t}^k-\hat{\bsx}_{i,t}^k) -{\boldsymbol{x}}_{avg,t}^{k} \right\| +\sumn \left\|  {\boldsymbol{x}}_{avg,t}^{k} - {\boldsymbol{x}}_{i,t}^{k} \right\|\nonumber \\
   %
   %
 &\leq  2 \sumn \left\|  {\boldsymbol{x}}_{i,t}^{k} -{\boldsymbol{x}}_{avg,t}^{k} \right\| +\alpha_t n M.
\end{flalign}
This, together with (\ref{proof_grad diff2 8_3}), implies that
\begin{flalign}
&\quad n\Gamma_1\sum\limits_{k=2}^{K_t} \sum\limits_{l=2}^k \sum\limits_{j=1}^n {\sigma_1}^{k-l} \left\|  \delta_{j,t}^{l-1} \right\| \nonumber \\
 & \leq  \frac{ 2n \Gamma_1 G_X}{1-{\sigma_1}} \sum\limits_{k=1}^{{K_t}-1} \sumn \left\|  {\boldsymbol{x}}_{i,t}^{k} -{\boldsymbol{x}}_{avg,t}^{k} \right\|+\frac{ n^2 \Gamma_1 G_X M}{1-{\sigma_1}} \alpha_t {K_t}.\label{proof_grad diff2 8_3-a1}
\end{flalign}
Then, substituting   (\ref{proof_grad diff2 8_1}),  (\ref{proof_grad diff2 8_2}) and (\ref{proof_grad diff2 8_3-a1}) into  (\ref{proof_grad diff2 7}), we have
\begin{flalign} \label{proof_grad diff2 10}
 &\quad \sumT \alpha_t\sum_{k=1}^{K_t} \sumn \left\|  \widehat{\nabla} f_{i,t}^k - \frac{1}{n}\sumn \nabla f_{i,t} (\hat{\bsx}_{i,t}^k)\right\| \nonumber \\
  &\leq \frac{ n^2 \Gamma_1 L_X}{1-{\sigma_1}} \sumT \alpha_t+ \frac{ n^2 \Gamma_1 G_X M}{1-{\sigma_1}}\sumT \alpha_t^2  {K_t}\nonumber \\
 &\quad + \frac{ 2n \Gamma_1 G_X}{1-{\sigma_1}} \sumT \sum\limits_{k=1}^{{K_t}-1} \sumn \alpha_t\left\|  {\boldsymbol{x}}_{i,t}^{k} -{\boldsymbol{x}}_{avg,t}^{k} \right\|.
\end{flalign}

On the other hand, similar to Lemma \ref{consistency2}, we easily obtain the following results:
\begin{flalign}
{\bsx}_{i,t}^{k+1}
&= \sum_{j=1}^n [A_t^k]_{ij}   {\bsx}_{j,t}^1 +\alpha_t \sum\limits_{l=1}^{k-1} \sum\limits_{j=1}^n [A_t^{k-l}]_{ij} \left(\boldsymbol{v}_{j,t}^l-\hat{\bsx}_{j,t}^l \right) \nonumber \\
 &\quad+ \alpha_t \left(\boldsymbol{v}_{i,t}^k-\hat{\bsx}_{i,t}^k \right),\label{proof_grad diff2 11_1} \\
  \bsx_{avg,t}^{k+1}
&=\frac{1}{n}\sumn {\boldsymbol{x}}_{i,t}^1 +\frac{\alpha_t}{n} \sum\limits_{l=1}^{k-1} \sumn ({\boldsymbol{v}}_{i,t}^l-\hat{\boldsymbol{x}}_{i,t}^l) \nonumber \\
 &\quad+ \alpha_t ({\boldsymbol{v}}_{avg,t}^k-{\boldsymbol{x}}_{avg,t}^k). \label{proof_grad diff2 11_2}
\end{flalign}
Based on these two conditions, we further obtain that
\begin{flalign}
&\quad \left\| {\bsx}_{i,t}^{k+1} -  \bsx_{avg,t}^{k+1}   \right\|\nonumber \\
&\leq   \sum_{j=1}^n \left| [A_t^k]_{ij}   -\frac{1}{n}\right|  \left\|{\boldsymbol{x}}_{j,t}^1 \right\| \nonumber \\
 &\quad + \alpha_t \sum\limits_{l=1}^{k-1} \sum\limits_{j=1}^n \left|   [A_t^{k-l}]_{ij} -  \frac{1}{n} \right\| \left\|{\boldsymbol{v}}_{i,t}^l-\hat{\boldsymbol{x}}_{i,t}^l\right\|+2\alpha_t M \nonumber \\
&\leq   \Gamma_1 {\sigma_1}^{k-1}\sum_{j=1}^n \left\|{\boldsymbol{x}}_{j,t}^1 \right\| +\alpha_t nM \Gamma_1\sum\limits_{l=1}^{k-1} {\sigma_1}^{k-l-1} +2\alpha_t M \nonumber \\
&=  \Gamma_1 {\sigma_1}^{k-1}\sum_{j=1}^n \left\|{\boldsymbol{x}}_{j,t}^1-{\boldsymbol{x}}_{j,1}+{\boldsymbol{x}}_{j,1}\right\| +\alpha_t n \Gamma_1 M \sum\limits_{s=1}^{k-1} {\sigma_1}^{s-1} \nonumber \\
 &\quad +2\alpha_t M \nonumber \\
&\leq nM\Gamma_1 {\sigma_1}^{k-1}+ \Gamma_1 {\sigma_1}^{k-1}\sum_{j=1}^n \left\|{\boldsymbol{x}}_{j,1}\right\| +\alpha_t n \Gamma_1 M \sum\limits_{s=1}^{k-1} {\sigma_1}^{s-1} \nonumber \\
 &\quad +2\alpha_t M. \label{proof_grad diff2 11_3}
\end{flalign}

Summing   (\ref{proof_grad diff2 11_3}) with respect to $i, k$ and $t$, we have
\begin{flalign} \label{proof_grad diff2 12}
&\quad \sum_{t=1}^T \sum_{k=1}^{K_t} \sum_{i=1}^n \alpha_t \|  {\bsx}_{i,t}^k -  {\bsx}_{avg,t}^k \| \nonumber \\
&\leq \sum_{t=1}^T\sum_{i=1}^n\alpha_t \|  {\bsx}_{i,t}^1 -  {\bsx}_{avg,t}^1 \|+\sum_{t=1}^T \sum_{k=2}^{K_t} \sum_{i=1}^n \alpha_t\|  {\bsx}_{i,t}^k -  {\bsx}_{avg,t}^k \|\nonumber \\
&\leq nM \sum_{t=1}^T \alpha_t +n^2M\Gamma_1 \sum_{t=1}^T \sum_{k=2}^{K_t}{\sigma_1}^{k-2} \alpha_t
\nonumber \\
 &\quad+ n\Gamma_1 \sum_{t=1}^T\sum_{k=2}^{K_t} \sum_{j=1}^n {\sigma_1}^{k-2} \alpha_t \| \boldsymbol{x}_{j,1}\| \nonumber \\
& \quad +  n^2 \Gamma_1 M  \sum_{t=1}^T \sum_{k=2}^{K_t}  \sum\limits_{s=1}^{k-1} {\sigma_1}^{s-1}\alpha_t^2  +2 n M  \sum_{t=1}^T \alpha_t^2 {K_t}\nonumber \\
&\leq  \left[nM+ \frac{n\Gamma_1 }{1-{\sigma_1}} \sum_{j=1}^n (\| \boldsymbol{x}_{j,1}\|+M) \right] \sum_{t=1}^T \alpha_t\nonumber \\
 &\quad +\left(\frac{n^2 \Gamma_1 M}{1-{\sigma_1}} + 2 n M\right) \sumT \alpha_t^2  {K_t}.
\end{flalign}

Combining   (\ref{proof_grad diff2 2}), (\ref{proof_grad diff2 10}) and (\ref{proof_grad diff2 12}) yields (\ref{grad diffience2-con1}).
%
The proof is complete.
\hfill$\square$


\section{Proof of Lemma \ref{lem F_k}}
Based on Algorithm \ref{algorithm 2} and the smooth property in Assumption \ref{assump: lips grad}, we have
\begin{flalign} \label{Proof lem F_k- 1}
 &\quad F_{t} \left(\boldsymbol{x}_{avg,t}^{k+1}\right) -F_{t}  \left(\boldsymbol{x}_{avg,t}^k \right) \nonumber \\
&\leq  \left< \nabla F_{t}  \left(\boldsymbol{x}_{avg,t}^k\right), \boldsymbol{x}_{avg,t}^{k+1}- \boldsymbol{x}_{avg,t}^k \right>  + \frac{nG_X}{2} \left\| \boldsymbol{x}_{avg,t}^{k+1}- \boldsymbol{x}_{avg,t}^k \right\|^2 \nonumber \\
&= \alpha_t \left< \nabla F_{t}  \left(\boldsymbol{x}_{avg,t}^k\right), \boldsymbol{v}_{avg,t}^{k}- \boldsymbol{x}_{avg,t}^k \right>  \nonumber \\
 &\quad + \frac{n\alpha_t^2 G_X}{2} \left\| \boldsymbol{v}_{avg,t}^{k}- \boldsymbol{x}_{avg,t}^k \right\|^2 \nonumber \\
&\leq\alpha_t \sumn \left< \frac{1}{n}\nabla F_{t}  \left(\boldsymbol{x}_{avg,t}^k\right), \boldsymbol{v}_{i,t}^{k}- \boldsymbol{x}_{avg,t}^k \right>  + \frac{n\alpha_t^2 G_XM^2}{2}
\end{flalign}
where we used the fact $\bsx_{avg,t}^{k+1}=\bsx_{avg,t}^{k}+\alpha_t (\bsv_{avg,t}^{k}-\bsx_{avg,t}^{k})$. Furthermore, it is noted that
\begin{flalign} \label{Proof lem F_k- 2}
& \quad \left< \frac{1}{n}\nabla F_{t}  \left(\boldsymbol{x}_{avg,t}^k\right), \boldsymbol{v}_{i,t}^{k}- \boldsymbol{x}_{avg,t}^k \right> \nonumber \\
&\leq  \left< \frac{1}{n}\nabla F_{t} \left(\boldsymbol{x}_{avg,t}^k\right)- \widehat{\nabla} f_{i,t}^k, \boldsymbol{v}_{i,t}^k- \boldsymbol{x}_{avg,t}^k \right> \nonumber \\
 &\quad+ \left< \widehat{\nabla} f_{i,t}^k, \boldsymbol{x}_t^*- \boldsymbol{x}_{avg,t}^k \right> \nonumber \\
&\leq  2M \left\| \frac{1}{n}\nabla F_{t} \left(\boldsymbol{x}_{avg,t}^k\right)- \widehat{\nabla} f_{i,t}^k \right\| \nonumber \\
 &\quad+ \left<\frac{1}{n}\nabla F_{t} \left(\boldsymbol{x}_{avg,t}^k\right), \boldsymbol{x}_t^*- \boldsymbol{x}_{avg,t}^k \right> \nonumber \\
&\leq  2M \left\| \frac{1}{n}\nabla F_{t} \left(\boldsymbol{x}_{avg,t}^k\right)- \widehat{\nabla} f_{i,t}^k \right\| + \frac{1}{n}\left[ F_{t} \left(\boldsymbol{x}_{t}^*\right)-F_{t} \left(\bsx_{avg,t}^k\right) \right]
\end{flalign}
where the first inequality and the last inequality utilize the optimality condition of the variable $\boldsymbol{v}_{i,t}^k$
and the convexity condition of $F_t( \bsx)$, respectively. Subtracting the term $ F_{t}(x_{t}^*)$ from both sides of   (\ref{Proof lem F_k- 1}) and substituting   (\ref{Proof lem F_k- 2}) into  (\ref{Proof lem F_k- 1}), the following inequality is obtained:
\begin{flalign} \label{Proof lem F_k- 3}
&\quad  F_{t} \left(\boldsymbol{x}_{avg,t}^{k+1}\right) - F_{t}(x_{t}^*) \nonumber \\
&= F_{t} \left(\boldsymbol{x}_{avg,t}^{k+1}\right) -F_{t}  \left(\boldsymbol{x}_{avg,t}^k \right)+F_{t}  \left(\boldsymbol{x}_{avg,t}^k \right)- F_{t}(x_{t}^*)
 \nonumber \\
&\leq  (1-\alpha_t)  \left[ F_{t}  \left(\boldsymbol{x}_{avg,t}^k \right)-F_{t}(x_{t}^*)\right] \nonumber \\
 &\quad+ 2\alpha_t M \sumn \left\| \frac{1}{n}\nabla F_{t} \left(\boldsymbol{x}_{avg,t}^k\right)- \widehat{\nabla} f_{i,t}^k \right\|    + \frac{n\alpha_t^2 G_XM^2}{2}\nonumber \\
&= (1-\alpha_t)^k  \left[ F_{t}  \left(\boldsymbol{x}_{avg,t}^1 \right)-F_{t}(x_{t}^*)\right]\nonumber \\
 &\quad + 2\alpha_t M \sum_{l=0}^{k-1} (1-\alpha_t)^l  \sumn \left\| \frac{1}{n}\nabla F_{t} \left(\boldsymbol{x}_{avg,t}^{k-l}\right)- \widehat{\nabla} f_{i,t}^{k-l} \right\| \nonumber\\
   &\quad+ \frac{n\alpha_t^2 G_XM^2}{2} \sum_{l=0}^{k-1}(1-\alpha_t)^l \nonumber \\
&\leq  (1-\alpha_t)^k  \left[ F_{t}  \left(\boldsymbol{x}_{avg,t}^1 \right)-F_{t}(x_{t}^*)\right] \nonumber \\
 &\quad+ 2\alpha_t M \sum_{l=1}^{k}   \sumn \left\| \frac{1}{n}\nabla F_{t} \left(\boldsymbol{x}_{avg,t}^{l}\right)- \widehat{\nabla} f_{i,t}^{l} \right\|    + \frac{n\alpha_t G_XM^2}{2}
\end{flalign}
where in the last inequality we utilize the two facts $(1-\alpha_t)^l \leq 1$ and $\sum_{l=0}^{k-1}(1-\alpha_t)^l\leq \frac{1}{\alpha_t}$. The proof is complete.
\hfill$\square$

\end{appendices}
\bibliographystyle{IEEEtran}
\bibliography{Arxiv_DOMFW_V2}

\begin{thebibliography}{10}
\providecommand{\url}[1]{#1}
\csname url@samestyle\endcsname
\providecommand{\newblock}{\relax}
\providecommand{\bibinfo}[2]{#2}
\providecommand{\BIBentrySTDinterwordspacing}{\spaceskip=0pt\relax}
\providecommand{\BIBentryALTinterwordstretchfactor}{4}
\providecommand{\BIBentryALTinterwordspacing}{\spaceskip=\fontdimen2\font plus
\BIBentryALTinterwordstretchfactor\fontdimen3\font minus
  \fontdimen4\font\relax}
\providecommand{\BIBforeignlanguage}[2]{{%
\expandafter\ifx\csname l@#1\endcsname\relax
\typeout{** WARNING: IEEEtran.bst: No hyphenation pattern has been}%
\typeout{** loaded for the language `#1'. Using the pattern for}%
\typeout{** the default language instead.}%
\else
\language=\csname l@#1\endcsname
\fi
#2}}
\providecommand{\BIBdecl}{\relax}
\BIBdecl

\bibitem{shalev2011online}
S.~Shalev-Shwartz, ``Online learning and online convex optimization,''
  \emph{Foundations and Trends in Machine Learning}, vol.~4, no.~2, pp.
  107--194, 2011.

\bibitem{hazan2016introduction}
E.~Hazan, ``Introduction to online convex optimization,'' \emph{Foundations and
  Trends{\textregistered} in Optimization}, vol.~2, no. 3-4, pp. 157--325,
  2016.

\bibitem{nedic2008distributed}
A.~Nedi{\'c}, A.~Olshevsky, A.~Ozdaglar, and J.~N. Tsitsiklis, ``Distributed
  subgradient methods and quantization effects,'' in \emph{2008 47th IEEE
  Conference on Decision and Control}, 2008, pp. 4177--4184.

\bibitem{nedic2018distributed}
A.~Nedi{\'c} and J.~Liu, ``Distributed optimization for control,'' \emph{Annual
  Review of Control, Robotics, and Autonomous Systems}, vol.~1, pp. 77--103,
  2018.

\bibitem{yang2019survey}
T.~Yang, X.~Yi, J.~Wu, Y.~Yuan, D.~Wu, Z.~Meng, Y.~Hong, H.~Wang, Z.~Lin, and
  K.~H. Johansson, ``A survey of distributed optimization,'' \emph{Annual
  Reviews in Control}, vol.~47, pp. 278--305, 2019.

\bibitem{li2022survey}
X.~Li, L.~Xie, and N.~Li, ``A survey of decentralized online learning,''
  \emph{arXiv preprint arXiv:2205.00473}, 2022.

\bibitem{10035518}
J.~Zhang, K.~You, and L.~Xie, ``Innovation compression for
  communication-efficient distributed optimization with linear convergence,''
  \emph{IEEE Transactions on Automatic Control}, 2023,
  doi:10.1109/TAC.2023.3241771.

\bibitem{liu2020unitary}
C.~Liu, H.~Li, and Y.~Shi, ``A unitary distributed subgradient method for
  multi-agent optimization with different coupling sources,''
  \emph{Automatica}, vol. 114, p. 108834, 2020.

\bibitem{li2021distributed}
W.~Li, X.~Zeng, Y.~Hong, and H.~Ji, ``Distributed consensus-based solver for
  semi-definite programming: An optimization viewpoint,'' \emph{Automatica},
  vol. 131, p. 109737, 2021.

\bibitem{xu2016distributed}
J.-M. Xu and Y.~C. Soh, ``A distributed simultaneous perturbation approach for
  large-scale dynamic optimization problems,'' \emph{Automatica}, vol.~72, pp.
  194--204, 2016.

\bibitem{6311406}
F.~Yan, S.~Sundaram, S.~Vishwanathan, and Y.~Qi, ``Distributed autonomous
  online learning: Regrets and intrinsic privacy-preserving properties,''
  \emph{IEEE Transactions on Knowledge and Data Engineering}, vol.~25, no.~11,
  pp. 2483--2493, 2013.

\bibitem{hosseini2013online}
S.~Hosseini, A.~Chapman, and M.~Mesbahi, ``Online distributed optimization via
  dual averaging,'' in \emph{52nd IEEE Conference on Decision and
  Control}.\hskip 1em plus 0.5em minus 0.4em\relax IEEE, 2013, pp. 1484--1489.

\bibitem{yuan2022distributedauto}
D.~Yuan, B.~Zhang, D.~W. Ho, W.~X. Zheng, and S.~Xu, ``Distributed online
  bandit optimization under random quantization,'' \emph{Automatica}, vol. 146,
  p. 110590, 2022.

\bibitem{yi2020distributedtsp}
X.~Yi, X.~Li, L.~Xie, and K.~H. Johansson, ``Distributed online convex
  optimization with time-varying coupled inequality constraints,'' \emph{IEEE
  Transactions on Signal Processing}, vol.~68, pp. 731--746, 2020.

\bibitem{9806334}
X.~Cao and T.~Ba\c{s}ar, ``Distributed constrained online convex optimization
  over multiple access fading channels,'' \emph{IEEE Transactions on Signal
  Processing}, vol.~70, pp. 3468--3483, 2022.

\bibitem{shahrampour2017distributed}
S.~Shahrampour and A.~Jadbabaie, ``Distributed online optimization in dynamic
  environments using mirror descent,'' \emph{IEEE Transactions on Automatic
  Control}, vol.~63, no.~3, pp. 714--725, 2017.

\bibitem{yuan2020distributed}
D.~Yuan, Y.~Hong, D.~W. Ho, and S.~Xu, ``Distributed mirror descent for online
  composite optimization,'' \emph{IEEE Transactions on Automatic Control},
  vol.~66, no.~2, pp. 714--729, 2020.

\bibitem{wang2020push}
C.~Wang, S.~Xu, D.~Yuan, B.~Zhang, and Z.~Zhang, ``Push-sum distributed online
  optimization with bandit feedback,'' \emph{IEEE Transactions on Cybernetics},
  vol.~52, no.~4, pp. 2263--2273, 2020.

\bibitem{zinkevich2003online}
M.~Zinkevich, ``Online convex programming and generalized infinitesimal
  gradient ascent,'' in \emph{Proceedings of the 20th international conference
  on machine learning (icml-03)}, 2003, pp. 928--936.

\bibitem{besbes2015non}
O.~Besbes, Y.~Gur, and A.~Zeevi, ``Non-stationary stochastic optimization,''
  \emph{Operations Research}, vol.~63, no.~5, pp. 1227--1244, 2015.

\bibitem{xu2022online}
Z.~Xu, H.~Zhou, and V.~Tzoumas, ``Online submodular coordination with bounded
  tracking regret: Theory, algorithm, and applications to multi-robot
  coordination,'' \emph{IEEE Robotics and Automation Letters}, vol.~8, no.~4,
  pp. 2261--2268, 2023.

\bibitem{pmlr-v70-zhang17g}
W.~Zhang, P.~Zhao, W.~Zhu, S.~C.~H. Hoi, and T.~Zhang, ``Projection-free
  distributed online learning in networks,'' in \emph{Proceedings of the 34th
  International Conference on Machine Learning}, 2017, pp. 4054--4062.

\bibitem{pmlr-v119-wan20b}
Y.~Wan, W.-W. Tu, and L.~Zhang, ``Projection-free distributed online convex
  optimization with ${O(\sqrt{T})}$ communication complexity,'' in
  \emph{Proceedings of the 37th International Conference on Machine Learning},
  2020, pp. 9818--9828.

\bibitem{wan2021projection}
Y.~Wan, G.~Wang, and L.~Zhang, ``Projection-free distributed online learning
  with strongly convex losses,'' \emph{arXiv preprint arXiv:2103.11102}, 2021.

\bibitem{thuang2022stochastic}
N.~K. Thang, A.~Srivastav, D.~Trystram, and P.~Youssef, ``A stochastic
  conditional gradient algorithm for decentralized online convex
  optimization,'' \emph{Journal of Parallel and Distributed Computing}, vol.
  169, pp. 334--351, 2022.

\bibitem{kalhan2021dynamic}
D.~S. Kalhan, A.~S. Bedi, A.~Koppel, K.~Rajawat, H.~Hassani, A.~K. Gupta, and
  A.~Banerjee, ``Dynamic online learning via {Frank-Wolfe} algorithm,''
  \emph{IEEE Transactions on Signal Processing}, vol.~69, pp. 932--947, 2021.

\bibitem{Wan_Xue_Zhang_2021}
Y.~Wan, B.~Xue, and L.~Zhang, ``Projection-free online learning in dynamic
  environments,'' in \emph{Proceedings of the AAAI Conference on Artificial
  Intelligence}, 2021, pp. 10\,067--10\,075.

\bibitem{zhang2023dynamic}
W.~Zhang, Y.~Shi, B.~Zhang, and D.~Yuan, ``Dynamic regret of distributed online
  {Frank-Wolfe} convex optimization,'' \emph{arXiv preprint arXiv:2302.00663},
  2023.

\bibitem{hazan2012projection}
E.~Hazan and S.~Kale, ``Projection-free online learning,'' in \emph{Proceedings
  of the 29th International Coference on International Conference on Machine
  Learning}, 2012, pp. 1843--1850.

\bibitem{hazan2008sparse}
E.~Hazan, ``Sparse approximate solutions to semidefinite programs,'' in
  \emph{Latin American Symposium on Theoretical Informatics}, 2008, pp.
  306--316.

\bibitem{harchaoui2015conditional}
Z.~Harchaoui, A.~Juditsky, and A.~Nemirovski, ``Conditional gradient algorithms
  for norm-regularized smooth convex optimization,'' \emph{Mathematical
  Programming}, vol. 152, no. 1-2, pp. 75--112, 2015.

\bibitem{7883821}
H.-T. Wai, J.~Lafond, A.~Scaglione, and E.~Moulines, ``Decentralized
  {Frank-Wolfe} algorithm for convex and nonconvex problems,'' \emph{IEEE
  Transactions on Automatic Control}, vol.~62, no.~11, pp. 5522--5537, 2017.

\bibitem{locatello2017unified}
F.~Locatello, R.~Khanna, M.~Tschannen, and M.~Jaggi, ``A unified optimization
  view on generalized matching pursuit and frank-wolfe,'' in \emph{Artificial
  Intelligence and Statistics}.\hskip 1em plus 0.5em minus 0.4em\relax PMLR,
  2017, pp. 860--868.

\bibitem{wu1983conditional}
Z.~Wu and K.~Teo, ``A conditional gradient method for an optimal control
  problem involving a class of nonlinear second-order hyperbolic partial
  differential equations,'' \emph{Journal of Mathematical Analysis and
  Applications}, vol.~91, no.~2, pp. 376--393, 1983.

\bibitem{zhang2017improved}
L.~Zhang, T.~Yang, J.~Yi, R.~Jin, and Z.-H. Zhou, ``Improved dynamic regret for
  non-degenerate functions,'' in \emph{Proceedings of the 31st International
  Conference on Neural Information Processing Systems}, 2017, pp. 732--741.

\bibitem{eshraghi2022dynamic}
N.~Eshraghi and B.~Liang, ``Dynamic regret bounds without lipschitz continuity:
  Online convex optimization with multiple mirror descent steps,'' in
  \emph{2022 American Control Conference (ACC)}.\hskip 1em plus 0.5em minus
  0.4em\relax IEEE, 2022, pp. 228--235.

\bibitem{wan2023improved}
Y.~Wan, L.~Zhang, and M.~Song, ``Improved dynamic regret for online
  frank-wolfe,'' \emph{arXiv preprint arXiv:2302.05620}, 2023.

\end{thebibliography}

\end{document}